\documentclass[9pt]{amsart}
\usepackage{amsmath}
\usepackage{amsxtra}
\usepackage{amscd}
\usepackage{amsthm}
\usepackage{amsfonts}
\usepackage{amssymb}
\usepackage{eucal}
\usepackage[all]{xy}
\usepackage{graphicx}

\newtheorem{cor}[subsubsection]{Corollary}
\newtheorem{lem}[subsubsection]{Lemma}
\newtheorem{prop}[subsubsection]{Proposition}

\newtheorem{conj}[subsubsection]{Conjecture}
\newtheorem{thm}[subsubsection]{Theorem}

\newtheorem{corconj}[subsubsection]{Corollary-of-Conjecture}

\theoremstyle{definition}

\theoremstyle{remark}
\newtheorem{rem}[subsubsection]{Remark}

\newcommand{\thmref}[1]{Theorem~\ref{#1}}

\newcommand{\secref}[1]{Sect.~\ref{#1}}

\newcommand{\propref}[1]{Proposition~\ref{#1}}
\newcommand{\corref}[1]{Corollary~\ref{#1}}
\newcommand{\conjref}[1]{Conjecture~\ref{#1}}

\numberwithin{equation}{section}

\newcommand{\nc}{\newcommand}
\nc{\renc}{\renewcommand}
\nc{\ssec}{\subsection}
\nc{\sssec}{\subsubsection}
\nc{\on}{\operatorname}

\nc\ol{\overline}
\nc\wt{\widetilde}
\nc\tboxtimes{\wt{\boxtimes}}
\nc\tstar{\wt{\star}}
\nc{\alp}{\alpha}

\nc{\ZZ}{{\mathbb Z}}
\nc{\NN}{{\mathbb N}}
\nc{\OO}{{\mathbb O}}
\renc{\SS}{{\mathbb S}}
\nc{\DD}{{\mathbb D}}
\nc{\GG}{{\mathbb G}}

\nc{\Fq}{{\mathbb F}_q}
\nc{\Fqb}{\ol{{\mathbb F}_q}}
\nc{\Ql}{\ol{{\mathbb Q}_\ell}}
\nc{\id}{\text{id}}
\nc\X{\mathcal X}

\nc{\Hom}{\on{Hom}}
\nc{\Lie}{\on{Lie}}
\nc{\Loc}{\on{Loc}}
\nc{\Pic}{\on{Pic}}
\nc{\Bun}{\on{Bun}}
\nc{\IC}{\on{IC}}
\nc{\Aut}{\on{Aut}}
\nc{\rk}{\on{rk}}
\nc{\Sh}{\on{Sh}}
\nc{\Perv}{\on{Perv}}
\nc{\pos}{{\on{pos}}}
\nc{\Conv}{\on{Conv}}
\nc{\Sph}{\on{Sph}}
\nc{\Sym}{\on{Sym}}
\nc{\BunBb}{\overline{\Bun}_B}
\nc{\BunNb}{\overline{\Bun}_N}
\nc{\BunTb}{\overline{\Bun}_T}
\nc{\BunBbm}{\overline{\Bun}_{B^-}}
\nc{\BunBbel}{\overline{\Bun}_{B,el}}
\nc{\BunBbmel}{\overline{\Bun}_{B^-,el}}
\nc{\Buno}{\overset{o}{\Bun}}
\nc{\BunPb}{{\overline{\Bun}_P}}
\nc{\BunBM}{\Bun_{B(M)}}
\nc{\BunBMb}{\overline{\Bun}_{B(M)}}
\nc{\BunPbw}{{\widetilde{\Bun}_P}}
\nc{\BunBP}{\widetilde{\Bun}_{B,P}}
\nc{\GUb}{\overline{G/U}}
\nc{\GUPb}{\overline{G/U(P)}}

\nc{\Hhom}{\underline{\on{Hom}}}
\nc\syminfty{\on{Sym}^{\infty}}
\nc\lal{\ol{\lambda}}
\nc\xl{\ol{x}}
\nc\thl{\ol{\theta}}
\nc\nul{\ol{\nu}}
\nc\mul{\ol{\mu}}
\nc\Sum\Sigma
\nc{\oX}{\overset{o}{X}{}}
\nc{\hl}{\overset{\leftarrow}h{}}
\nc{\hr}{\overset{\rightarrow}h{}}
\nc{\M}{{\mathcal M}}
\nc{\N}{{\mathcal N}}
\nc{\F}{{\mathcal F}}
\nc{\D}{{\mathcal D}}
\nc{\Q}{{\mathcal Q}}
\nc{\Y}{{\mathcal Y}}
\nc{\G}{{\mathcal G}}
\nc{\E}{{\mathcal E}}
\nc{\CalC}{{\mathcal C}}
\nc\Dh{\widehat{\D}}

\nc{\C}{{\mathcal C}}
\nc{\K}{{\mathcal K}}
\renewcommand{\H}{{\mathcal H}}

\nc{\T}{{\mathcal T}}
\nc{\V}{{\mathcal V}}
\renc{\P}{{\mathcal P}}
\nc{\A}{{\mathcal A}}
\nc{\B}{{\mathcal B}}
\nc{\U}{{\mathcal U}}

\nc{\Gr}{{\on{Gr}}}

\nc{\frn}{{\check{\mathfrak u}(P)}}

\nc{\fC}{\mathfrak C}
\nc{\p}{\mathfrak p}
\nc{\q}{\mathfrak q}
\nc\f{{\mathfrak f}}

\nc{\qo}{{\mathfrak q}}
\nc{\po}{{\mathfrak p}}
\nc{\s}{{\mathfrak s}}
\nc\w{\text{w}}

\renewcommand{\mod}{{\on{-mod}}}
\newcommand{\comod}{{\on{-comod}}}

\nc\mathi\iota
\nc\Spec{\on{Spec}}
\nc\Mod{\on{Mod}}
\nc{\tw}{\widetilde{\mathfrak t}}
\nc{\pw}{\widetilde{\mathfrak p}}
\nc{\qw}{\widetilde{\mathfrak q}}
\nc{\jw}{\widetilde j}

\nc{\grb}{\overline{\Gr}}
\nc{\I}{\mathcal I}

\nc{\lambdach}{{\check\lambda}}
\nc{\Lambdach}{{\check\Lambda}{}}
\nc{\much}{{\check\mu}}
\nc{\omegach}{{\check\omega}}
\nc{\nuch}{{\check\nu}}
\nc{\etach}{{\check\eta}}
\nc{\alphach}{{\check\alpha}}
\nc{\oblvtach}{{\check\oblvta}}
\nc{\rhoch}{{\check\rho}}
\nc{\ch}{{\check h}}

\nc{\Hb}{\overline{\H}}

\nc{\BA}{{\mathbb{A}}}
\nc{\BC}{{\mathbb{C}}}
\nc{\BE}{{\mathbb{E}}}
\nc{\BF}{{\mathbb{F}}}
\nc{\BG}{{\mathbb{G}}}
\nc{\BM}{{\mathbb{M}}}
\nc{\BO}{{\mathbb{O}}}
\nc{\BD}{{\mathbb{D}}}
\nc{\BL}{{\mathbb{L}}}
\nc{\Bl}{{\mathbb{l}}}
\nc{\BN}{{\mathbb{N}}}
\nc{\BP}{{\mathbb{P}}}
\nc{\BQ}{{\mathbb{Q}}}
\nc{\BR}{{\mathbb{R}}}
\nc{\BZ}{{\mathbb{Z}}}
\nc{\BS}{{\mathbb{S}}}

\nc{\CA}{{\mathcal{A}}}
\nc{\CB}{{\mathcal{B}}}

\nc{\CE}{{\mathcal{E}}}
\nc{\CF}{{\mathcal{F}}}
\nc{\CH}{{\mathcal{H}}}

\nc{\CL}{{\mathcal{L}}}
\nc{\CC}{{\mathcal{C}}}
\nc{\CG}{{\mathcal{G}}}
\nc{\CM}{{\mathcal{M}}}
\nc{\CN}{{\mathcal{N}}}
\nc{\CK}{{\mathcal{K}}}
\nc{\CO}{{\mathcal{O}}}
\nc{\CP}{{\mathcal{P}}}
\nc{\CQ}{{\mathcal{Q}}}
\nc{\CR}{{\mathcal{R}}}
\nc{\CS}{{\mathcal{S}}}
\nc{\CT}{{\mathcal{T}}}
\nc{\CU}{{\mathcal{U}}}
\nc{\CV}{{\mathcal{V}}}
\nc{\CW}{{\mathcal{W}}}
\nc{\CX}{{\mathcal{X}}}
\nc{\CY}{{\mathcal{Y}}}
\nc{\CZ}{{\mathcal{Z}}}
\nc{\CI}{{\mathcal{I}}}
\nc{\cD}{{\mathcal{D}}}
\nc{\ocD}{\overset{\circ}{\mathcal{D}}}

\nc{\csM}{{\check{\mathcal A}}{}}
\nc{\oM}{{\overset{\circ}{\mathcal M}}{}}
\nc{\obM}{{\overset{\circ}{\mathbf M}}{}}
\nc{\oCA}{{\overset{\circ}{\mathcal A}}{}}
\nc{\obA}{{\overset{\circ}{\mathbf A}}{}}
\nc{\ooM}{{\overset{\circ}{M}}{}}
\nc{\osM}{{\overset{\circ}{\mathsf M}}{}}
\nc{\vM}{{\overset{\bullet}{\mathcal M}}{}}
\nc{\nM}{{\underset{\bullet}{\mathcal M}}{}}
\nc{\oD}{{\overset{\circ}{\mathcal D}}{}}
\nc{\obD}{{\overset{\circ}{\mathbf D}}{}}
\nc{\oA}{{\overset{\circ}{\mathbb A}}{}}
\nc{\op}{{\overset{\bullet}{\mathbf p}}{}}
\nc{\cp}{{\overset{\circ}{\mathbf p}}{}}
\nc{\oU}{{\overset{\bullet}{\mathcal U}}{}}
\nc{\oZ}{{\overset{\circ}{\mathcal Z}}{}}
\nc{\ofZ}{{\overset{\circ}{\mathfrak Z}}{}}
\nc{\oF}{{\overset{\circ}{\fF}}}

\nc{\fa}{{\mathfrak{a}}}
\nc{\fb}{{\mathfrak{b}}}
\nc{\fc}{{\mathfrak{c}}}
\nc{\fch}{{\mathfrak{ch}}}
\nc{\fd}{{\mathfrak{d}}}
\nc{\ff}{{\mathfrak{f}}}
\nc{\fg}{{\mathfrak{g}}}
\nc{\fgl}{{\mathfrak{gl}}}
\nc{\fh}{{\mathfrak{h}}}
\nc{\fj}{{\mathfrak{j}}}
\nc{\fl}{{\mathfrak{l}}}
\nc{\fm}{{\mathfrak{m}}}
\nc{\fn}{{\mathfrak{n}}}
\nc{\fu}{{\mathfrak{u}}}
\nc{\fp}{{\mathfrak{p}}}
\nc{\fr}{{\mathfrak{r}}}
\nc{\fs}{{\mathfrak{s}}}
\nc{\ft}{{\mathfrak{t}}}
\nc{\fT}{{\mathfrak{T}}}
\nc{\fz}{{\mathfrak{z}}}
\nc{\fsl}{{\mathfrak{sl}}}
\nc{\hsl}{{\widehat{\mathfrak{sl}}}}
\nc{\hgl}{{\widehat{\mathfrak{gl}}}}
\nc{\hg}{{\widehat{\mathfrak{g}}}}
\nc{\htt}{{\widehat{\mathfrak{t}}}}
\nc{\chg}{{\widehat{\mathfrak{g}}}{}^\vee}
\nc{\hn}{{\widehat{\mathfrak{n}}}}
\nc{\chn}{{\widehat{\mathfrak{n}}}{}^\vee}

\nc{\fA}{{\mathfrak{A}}}
\nc{\fB}{{\mathfrak{B}}}
\nc{\fD}{{\mathfrak{D}}}
\nc{\fE}{{\mathfrak{E}}}
\nc{\fF}{{\mathfrak{F}}}
\nc{\fG}{{\mathfrak{G}}}
\nc{\fK}{{\mathfrak{K}}}
\nc{\fL}{{\mathfrak{L}}}
\nc{\fM}{{\mathfrak{M}}}
\nc{\fN}{{\mathfrak{N}}}
\nc{\fP}{{\mathfrak{P}}}
\nc{\fU}{{\mathfrak{U}}}
\nc{\fV}{{\mathfrak{V}}}
\nc{\fZ}{{\mathfrak{Z}}}

\nc{\bb}{{\mathbf{b}}}
\nc{\bc}{{\mathbf{c}}}
\nc{\bd}{{\mathbf{d}}}
\nc{\bbf}{{\mathbf{f}}}
\nc{\be}{{\mathbf{e}}}
\nc{\bg}{{\mathbf{g}}}
\nc{\bi}{{\mathbf{i}}}
\nc{\bj}{{\mathbf{j}}}
\nc{\bn}{{\mathbf{n}}}
\nc{\bp}{{\mathbf{p}}}
\nc{\bq}{{\mathbf{q}}}
\nc{\bu}{{\mathbf{u}}}
\nc{\bv}{{\mathbf{v}}}
\nc{\bx}{{\mathbf{x}}}
\nc{\bs}{{\mathbf{s}}}
\nc{\by}{{\mathbf{y}}}
\nc{\bw}{{\mathbf{w}}}
\nc{\bA}{{\mathbf{A}}}
\nc{\bK}{{\mathbf{K}}}
\nc{\bB}{{\mathbf{B}}}
\nc{\bC}{{\mathbf{C}}}
\nc{\bG}{{\mathbf{G}}}
\nc{\bD}{{\mathbf{D}}}
\nc{\bH}{{\mathbf{He}}}
\nc{\bM}{{\mathbf{M}}}
\nc{\bN}{{\mathbf{N}}}
\nc{\bO}{{\mathbf{O}}}
\nc{\bV}{{\mathbf{V}}}
\nc{\bW}{{\mathbf{Wh}}}
\nc{\bX}{{\mathbf{X}}}
\nc{\bY}{{\mathbf{Y}}}
\nc{\bZ}{{\mathbf{Z}}}
\nc{\bS}{{\mathbf{S}}}
\nc{\bT}{{\mathbf{T}}}

\nc{\sA}{{\mathsf{A}}}
\nc{\sB}{{\mathsf{B}}}
\nc{\sC}{{\mathsf{C}}}
\nc{\sD}{{\mathsf{D}}}
\nc{\sF}{{\mathsf{F}}}
\nc{\sG}{{\mathsf{G}}}
\nc{\sK}{{\mathsf{K}}}
\nc{\sM}{{\mathsf{M}}}
\nc{\sO}{{\mathsf{O}}}
\nc{\sU}{{\mathsf{U}}}
\nc{\sW}{{\mathsf{W}}}
\nc{\sQ}{{\mathsf{Q}}}
\nc{\sP}{{\mathsf{P}}}
\nc{\sY}{{\mathsf{Y}}}
\nc{\sZ}{{\mathsf{Z}}}
\nc{\sfp}{{\mathsf{p}}}
\nc{\sfq}{{\mathsf{q}}}
\nc{\sr}{{\mathsf{r}}}
\nc{\sk}{{\mathsf{k}}}
\nc{\su}{{\mathsf{u}}}
\nc{\sv}{{\mathsf{v}}}
\nc{\sg}{{\mathsf{g}}}
\nc{\sff}{{\mathsf{f}}}
\nc{\sfb}{{\mathsf{b}}}
\nc{\sfc}{{\mathsf{c}}}
\nc{\sd}{{\mathsf{d}}}

\nc{\BK}{{\bar{K}}}

\nc{\tA}{{\widetilde{\mathbf{A}}}}
\nc{\tB}{{\widetilde{\mathcal{B}}}}
\nc{\tg}{{\widetilde{\mathfrak{g}}}}
\nc{\tG}{{\widetilde{G}}}
\nc{\TM}{{\widetilde{\mathbb{M}}}{}}
\nc{\tO}{{\widetilde{\mathsf{O}}}{}}
\nc{\tU}{{\widetilde{\mathfrak{U}}}{}}
\nc{\TZ}{{\tilde{Z}}}
\nc{\tx}{{\tilde{x}}}
\nc{\tbv}{{\tilde{\bv}}}
\nc{\tfP}{{\widetilde{\mathfrak{P}}}{}}
\nc{\tz}{{\tilde{\zeta}}}
\nc{\tmu}{{\tilde{\mu}}}

\nc{\urho}{\underline{\rho}}
\nc{\uB}{\underline{B}}
\nc{\uC}{{\underline{\mathbb{C}}}}
\nc{\ui}{\underline{i}}
\nc{\uj}{\underline{j}}
\nc{\ofP}{{\overline{\mathfrak{P}}}}
\nc{\oB}{{\overline{\mathcal{B}}}}
\nc{\og}{{\overline{\mathfrak{g}}}}
\nc{\oI}{{\overline{I}}}

\nc{\eps}{\varepsilon}
\nc{\hrho}{{\hat{\rho}}}

\nc{\one}{{\mathbf{1}}}
\nc{\two}{{\mathbf{t}}}

\nc{\Rep}{{\mathop{\operatorname{\rm Rep}}}}
\nc{\Tot}{{\mathop{\operatorname{\rm Tot}}}}
\nc{\Ker}{{\mathop{\operatorname{\rm Ker}}}}
\nc{\Hilb}{{\mathop{\operatorname{\rm Hilb}}}}
\nc{\End}{{\mathop{\operatorname{\rm End}}}}
\nc{\Ext}{{\mathop{\operatorname{\rm Ext}}}}
\nc{\CHom}{{\mathop{\operatorname{{\mathcal{H}}\it om}}}}
\nc{\GL}{{\mathop{\operatorname{\rm GL}}}}
\nc{\gr}{{\mathop{\operatorname{\rm gr}}}}
\nc{\Id}{{\mathop{\operatorname{\rm Id}}}}
\nc{\de}{{\mathop{\operatorname{\rm def}}}}
\nc{\length}{{\mathop{\operatorname{\rm length}}}}
\nc{\supp}{{\mathop{\operatorname{\rm supp}}}}

\nc{\Cliff}{{\mathsf{Cliff}}}
\nc{\Fl}{\on{Fl}}
\nc{\Fib}{{\mathsf{Fib}}}
\nc{\Coh}{{\on{Coh}}}
\nc{\QCoh}{{\on{QCoh}}}
\nc{\IndCoh}{{\on{IndCoh}}}
\nc{\FCoh}{{\mathsf{FCoh}}}

\nc{\reg}{{\text{\rm reg}}}

\nc{\cplus}{{\mathbf{C}_+}}
\nc{\cminus}{{\mathbf{C}_-}}
\nc{\cthree}{{\mathbf{C}_*}}
\nc{\Qbar}{{\bar{Q}}}
\nc\Eis{{\on{Eis}}}
\nc\Eisb{\ol\Eis{}}
\nc\Eisr{\on{Eis}^{rat}{}}
\nc\wh{\widehat}
\nc{\Def}{\on{Def_{\check{\fb}}(E)}}
\nc{\barZ}{\overline{Z}{}}
\nc{\barbarZ}{\overline{\barZ}{}}
\nc{\barpi}{\overline\pi}
\nc{\barbarpi}{\overline\barpi}
\nc{\barpip}{\overline\pi{}^+}
\nc{\barpim}{\overline\pi{}^-}

\nc{\fq}{\mathfrak q}

\nc{\fqb}{\ol{\fq}{}}
\nc{\fpb}{\ol{\fp}{}}
\nc{\fpr}{{\fp^{rat}}{}}
\nc{\fqr}{{\fq^{rat}}{}}

\nc{\hattimes}{\wh\otimes}

\nc{\bh}{{{\mathbf h}}}
\nc{\bk}{{{\mathbf k}}}
\nc{\bOmega}{{\overline{\Omega(\check \fn)}}}

\nc{\seq}[1]{\stackrel{#1}{\sim}}

\nc{\cT}{{\check{T}}}
\nc{\cG}{{\check{G}}}
\nc{\cM}{{\check{M}}}
\nc{\cB}{{\check{B}}}
\nc{\cP}{{\check{P}}}

\nc{\ct}{{\check{\mathfrak t}}}
\nc{\cg}{{\check{\fg}}}
\nc{\cb}{{\check{\fb}}}
\nc{\cn}{{\check{\fn}}}

\nc{\cLambda}{{\check\Lambda}}

\nc{\cla}{{\check\lambda}}
\nc{\cmu}{{\check\mu}}
\nc{\cnu}{{\check\nu}}
\nc{\ceta}{{\check\eta}}

\nc{\DefbE}{{\on{Def}_{\cB}(E_\cT)}}

\nc{\imathb}{{\ol{\imath}}}
\nc{\rlr}{\overset{\longrightarrow}{\underset{\longrightarrow}\longleftarrow}}

\nc{\oBun}{\overset{\circ}\Bun}
\nc{\oSht}{\overset{\circ}\Sht}
\nc{\LocSys}{\on{LocSys}}
\nc{\BunBbb}{\ol{\ol{Bun}}_B}
\nc{\BunBr}{\Bun_B^{rat}}
\nc{\BunBrp}{\Bun_B^{rat,polar}}
\nc{\BunTrp}{\Bun_T^{rat,polar}}
\nc{\BunNr}{\Bun_N^{rat}}
\nc{\BunNre}{\Bun_N^{enh,rat}}
\nc{\BunTr}{\Bun_T^{rat}}
\nc{\Vect}{\on{Vect}}
\nc{\Whit}{\on{Whit}}
\nc{\CTb}{\ol{\on{CT}}}
\nc{\Ran}{\on{Ran}}
\nc{\CTr}{\on{CT}^{rat}{}}
\nc\jmathr{\jmath^{rat}{}}
\nc{\ux}{\underline{x}}
\nc{\clambda}{{\check\lambda}}
\nc{\calpha}{{\check\alpha}}
\nc{\ind}{{\mathbf{ind}}}
\nc{\oblv}{{\mathbf{oblv}}}
\nc{\coeff}{\on{W-coeff}}
\nc{\Poinc}{\on{Poinc}}
\nc{\Dmod}{\on{D-mod}}
\nc{\dr}{\on{dR}}
\nc{\oCZ}{\overset{\circ}\CZ}
\nc{\KL}{\on{KL}}
\nc{\triv}{{\mathbf{triv}}}

\nc{\dgSch}{\on{DGSch}}
\nc{\Sch}{\on{Sch}}
\nc{\affdgSch}{\on{DGSch}^{\on{aff}}}
\nc{\affSch}{\on{Sch}^{\on{aff}}}
\nc{\Sing}{\on{Sing}}
\nc{\inftygroup}{\infty\on{-Grpd}}
\renc{\dr}{{\on{dr}}}
\nc\Maps{\on{Maps}}
\nc\Res{\on{Res}}
\nc\bMaps{\mathbf{Maps}}
\nc{\ul}{\underline}
\nc{\bNP}{\mathbf{N(P)}}
\nc{\ofc}{\overset{\circ}\fch}
\nc{\ppart}{(\!(t)\!)}
\nc{\qqart}{[\![t]\!]}
\nc{\crit}{\on{crit}}
\nc{\DGCat}{\on{DGCat}}
\nc{\Shv}{\on{Shv}}
\nc{\bDelta}{\mathbf{\Delta}}
\nc{\genB}{{\overset{\on{gen}}\to B}}
\nc{\genP}{{\underset{\on{gen}}\longrightarrow P}}
\nc{\genN}{{\underset{\on{gen}}\longrightarrow N}}
\nc{\semiinf}{{\frac{\infty}{2}+\bullet}}
\nc{\mmod}{\on{-}\mathbf{mod}}
\nc{\AdFr}{\on{Ad}_{\on{Frob}}}
\nc{\Frob}{\on{Frob}}
\nc{\Tr}{\on{Tr}}
\nc{\Sht}{\on{Sht}}
\nc{\sfe}{\mathsf{e}}
\nc{\tCat}{{2\on{-Cat}}}
\nc{\tIndCoh}{{2\on{-IndCoh}}}
\nc{\tQCoh}{{2\on{-QCoh}}}
\nc{\uShv}{\underline{\Shv}}
\nc{\qLisse}{\on{QLisse}}
\nc{\Nilp}{{\on{Nilp}}}
\nc{\sotimes}{\overset{!}\otimes} 
\nc{\AGCat}{\on{AGCat}}
\nc{\LS}{\on{LS}}
\nc{\Cat}{\on{-Cat}}
\nc{\bo}{\mathbf o}

\begin{document}

%\begin{quote}
%{\small ``In jedem Minus steckt ein Plus. Vielleicht habe ich so etwas gesagt, aber man braucht das doch nicht allzu w\"ortlich zu nehmen."}
%\hskip1.3cm {\tiny R.~Musil. Der Mann ohne Eigenschaften.}
%\end{quote}

\vskip1cm

\title{Local and global Langlands conjecture(s)\\  over function fields}

\author{Dennis Gaitsgory} 

%\begin{abstract}
%This paper is a write-up of the plenary talk at ICM, 2026.
%\end{abstract}

\date{\today}

\maketitle

\centerline{\bf Write-up of the plenary talk at ICM, 2026}

\section*{Introduction}

\ssec{What is this paper about?}

The goal of this paper is to propose a set of conjectures whose aim is to answer the basic question of the 
Langlands program (over function fields): how to describe the space of automorphic functions in terms of the spectral side (i.e., Langlands parameters)? 

\smallskip

The approach we take can be summarized as the following three steps:

\begin{itemize}

\item Formulate (a higher-categorical) geometric statement;

\item Apply the operation of (higher-categorical) trace of Frobenius, and identify the result 
with the classical object of interest;

\item Deduce the sought-for description for the classical object.

\end{itemize}

\ssec{Structure}

The main body of this paper consists of three sections:

\smallskip

In \secref{s:unram} we review some recent progress in the global unramified theory, following the papers
\cite{AGKRRV1,AGKRRV2,AGKRRV3} and \cite{GR}. 

\smallskip

In \secref{s:local} we propose a set of conjectures pertaining to the local Langlands theory. 

\smallskip

In \secref{s:global} we propose a set of conjectures pertaining to the global Langlands theory
(whose formulation takes the local theory as an input). 

\ssec{Notation}

Throughout this paper, $G$ will be a reductive group and $X$ a smooth and complete curve over a ground field $k$,
assumed algebraically closed. We will denote by $\Shv(-)$ the sheaf theory of $\ol\BQ_\ell$-adic sheaves
(see \cite[Sect. 1.1]{AGKRRV1}). 

\smallskip

We will work with \emph{higher algebra} and \emph{derived algebraic geometry} over the field of coefficients $\sfe:=\ol\BQ_\ell$. 
In particular, we will work with $\sfe$-linear DG categories (see \cite[Chapter 1, Sect.10]{GaRo} for what we mean by this). We will
denote by $\Vect$ the DG category of $\sfe$-vector spaces. We refer the reader to \cite[Chapter 2]{GaRo} for our conventions
regarding algebraic geometry (e.g., the notion of prestack), and to Chapter 3 of {\it loc. cit.} for $\QCoh(-)$. 

\smallskip

In all other respects, our notations and conventions follow those of \cite{AGKRRV1}.

\ssec{Acknowledgements}

The contents of this paper describe work-in-progress joint with S.~Raskin; some parts are joint also with other mathematicians:
D.~Arinkin, G.~Dhillon, A.~Eteve, A.~Genestier, D.~Kazhdan, V.~Lafforgue, N.~Rozenblyum Y.~Varshavsky and D.~Yang.  

\smallskip

The author wishes also to thank A.~Beilinson, D.~Ben-Zvi, V.~Drinfeld, D.~Nadler, W.~Sawin, P.~Scholze, C.~Xue, Z.~Yun
and X.~Zhu, whose ideas have informed his thinking about the subject.

\section{The global unramified case} \label{s:unram}

In this section we review the global unramified geometric Langlands conjecture and 
its application to the theory of automorphic functions in the global unramified case
over function fields.

\ssec{The players} \label{ss:BunG Nilp}

\sssec{} \label{sss:goals}

Let $\Bun_G$ denote the moduli stack of $G$-bundles on $X$. In the geometric theory 
our main object of study is the category $\Shv(\Bun_G)$, the derived category of
$\ol\BQ_\ell$-adic sheaves on $\Bun_G$. 

\smallskip

Here are our two goals:

\begin{enumerate}

\item Relate $\Shv(\Bun_G)$ to the space of unramified automorphic functions (when $k=\ol\BF_q$,
but $G$ and $X$ are defined over $\BF_q$);

\item Relate $\Shv(\Bun_G)$ to the \emph{spectral side}, i.e.,
to the ``space" of $\cG$-local systems on $X$, where $\cG$ is the Langlands-dual of $G$,
which is a reductive group over $\ol\BQ_\ell=:\sfe$.

\end{enumerate}

If we could achieve these two goals, we could hope to relate the space of automorphic functions
to $\cG$-local systems, thereby realizing Langlands' vision.

\smallskip

Here by the space of (unramified, compactly supported, $\sfe$-valued) automorphic functions we mean the vector space
$$\on{Autom}(X,G):=\on{Funct}_c(\Bun_G(\BF_q)).$$

\sssec{}

There is an ``obvious" strategy of how one may tackle goal (1) above. Namely, since $G$ and $X$
are defined over $\BF_q$, so is $\Shv(\Bun_G)$, and hence it carries an endomorphism, given by
\emph{geometric Frobenius}, to be denoted $\Frob$. 

\smallskip

This endomorphism induces an endofunctor
(in fact, a self-equivalence) $\Frob_*$ on $\Shv(\Bun_G)$. and we can consider the vector space
$$\Tr(\Frob_*,\Shv(\Bun_G))\in \Vect,$$
see \secref{ss:cat Tr} for the reminder on the formalism of categorical trace. 

\smallskip

More abstractly, let $\CY$ be an algebraic stack over $\ol\BF_q$ but defined over $\BF_q$, and consider 
$$\Tr(\Frob_*,\Shv(\CY))\in \Vect.$$

If $\CY$ is quasi-compact, a construction from \cite[Sect. 22.2]{AGKRRV1} defines a map
\begin{equation} \label{e:Tr Frob on sheaves qc}
\on{LT}^{\on{true}}:\Tr(\Frob_*,\Shv(\CY))\to \on{Funct}(\CY(\BF_q)),
\end{equation}
where $\on{LT}$ stands for ``local term". 

\smallskip

When $\CY$ is not quasi-compact, instead of \eqref{e:Tr Frob on sheaves qc} we rather have a map
\begin{equation} \label{e:Tr Frob on sheaves}
\on{LT}^{\on{true}}:\Tr(\Frob_*,\Shv(\CY))\to \on{Funct}_c(\CY(\BF_q)),
\end{equation}
where superscript ``c" stands for functions with finite support. 

\sssec{}

Here is how \eqref{e:Tr Frob on sheaves} is related to Grothendieck's \emph{functions-sheaves} dictionary: let
$\CF$ be a compact (and in particular, constructible) object of $\Shv(\CY)$, equipped with a \emph{weak Weil structure},
i.e., a map
$$\CF\overset{\alpha}\to \Frob_*(\CF).$$

\smallskip

On the one hand, taking traces of the Frobenius on the stalks of $\CF$, we obtain a function
$$\on{funct}(\CF,\alpha)\in \on{Funct}_c(\CY(\BF_q)).$$

On the other hand, to the pair $(\CF,\alpha)$ one can associate its class 
$$\on{cl}(\CF,\alpha)\in \Tr(\Frob_*,\Shv(\CY)),$$
see \secref{sss:class}. The claim\footnote{Proved in \cite[Corollary 0.9]{GV}.} is that 
$\on{funct}(\CF,\alpha)=\on{LT}^{\on{true}}(\on{cl}(\CF,\alpha)$. 

\begin{rem}

In fact, the assignment $(\CF,\alpha)\mapsto \on{funct}(\CF,\alpha)$
can be upgraded to a map (see \cite[Sect. 22.1.4]{AGKRRV1})
$$\on{LT}^{\on{naive}}:\Tr(\Frob_*,\Shv(\CY))\to \on{Funct}_c(\CY(\BF_q)),$$
and the statement is that
$\on{LT}^{\on{naive}}=\on{LT}^{\on{true}}$
as maps 
$\Tr(\Frob_*,\Shv(\CY))\rightrightarrows \on{Funct}_c(\CY(\BF_q))$. 

\end{rem} 

\sssec{} \label{sss:fail goal}

However, the map $\on{LT}^{\on{true}}$ of \eqref{e:Tr Frob on sheaves} is \emph{not} an isomorphism\footnote{Unless $\CY$ is \emph{restricted}, see \secref{sss:restr stack}.},
unless we understand $\Tr(\Frob_*,-)$ in the sense of $\AGCat$, see \secref{sss:trace in AGCat}. So considering $\Tr(\Frob_*,\Shv(\Bun_G))$ does not quite 
help us achieve goal (1) from \secref{sss:goals}. 

\smallskip

Moreover, when we look at $\Shv(\Bun_G)$, it is not clear how to even start tackling goal (2). We do not know what the relevant algebro-geometric
object should be that can reasonably be thought of as the moduli stack of $\cG$-local systems. 

\smallskip 

A way to resolve both these issues has been suggested in the series of papers \cite{AGKRRV1}-\cite{AGKRRV3}. Namely, in {\it loc. cit.} it was suggested
that instead of the entire category $\Shv(\Bun_G)$, one should consider its full subcategory, denoted $\Shv_\Nilp(\Bun_G)$, singled out by the
singular support condition.

\sssec{}

Namely, we recall that following \cite{Bei}, given an algebraic stack $\CY$ and a constructible object $\CF\in \Shv(\CY)$,
one can assign to it a conical subset 
$$\on{SingSupp}(\CF)\subset T^*(\CY).$$

Vice versa, given a conical subset $\CN\subset T^*(\CY)$, one can consider the full subcategory
$$\Shv_\CN(\CY)\subset \Shv(\CY),$$
consisting of objects, all of whose perverse cohomologies have subquotients that lie in $\CN$.

\sssec{}

We take $\CY=\Bun_G$, and we identify $T^*(\Bun_G)$ with the moduli space of Higgs bundles. We let
$\Nilp\subset T^*(\Bun_G)$ be the global nilpotent cone, defined as in \cite[Sect. 14.1.2]{AGKRRV1}. We let 
\begin{equation} \label{e:Nilp as subcat}
\Shv_\Nilp(\Bun_G)\subset \Shv(\Bun_G)
\end{equation} 
denote the corresponding full subcategory. 

\begin{rem}

We emphasize that the relevance of the subcategory $\Shv_\Nilp(\Bun_G)\subset \Shv(\Bun_G)$ was realized early in the development 
of the subject, namely, already in G.~Laumon's paper \cite{Laum} (see Sect. 6 in {\it loc. cit.}). 

\end{rem}

\sssec{} \label{sss:Hecke lisse}

One the key theorems in the paper \cite{AGKRRV1}, namely, Theorem 14.4.3 in {\it loc. cit.}, says that the subcategory  
\eqref{e:Nilp as subcat} can also be characterized differently.

\smallskip

Namely, recall that the (naive) geometric Satake construction assigns to every $V\in \Rep(\cG)$ a \emph{Hecke functor}
\begin{equation} \label{e:Hecke functor}
\on{H}_V:\Shv(\Bun_G)\to \Shv(\Bun_G\times X).
\end{equation} 

\smallskip

Let\footnote{In the formula below, ``HL" stands for ``Hecke-lisse".} 
$$\Shv_{\on{HL}}(\Bun_G)\subset \Shv(\Bun_G)$$
denote the full subcategory consisting of objects $\CF$, such that for any $V\in \Rep(\cG)$, the resulting object 
$\on{H}_V(\CF)$ belongs to the full subcategory
$$\Shv(\Bun_G)\otimes \qLisse(X)\subset \Shv(\Bun_G)\otimes \Shv(X)\subset \Shv(\Bun_G\times X),$$
where 
$\qLisse(X) \subset \Shv(X)$ is the full subcategory consisting of objects whose perverse cohomology are 
(ind)\emph{lisse} (see \cite[Definition 1.2.6]{AGKRRV1}).

\smallskip

The result of \cite{AGKRRV1}\footnote{The inclusion $\subseteq$ had been earlier established in \cite{NY}.} 
mentioned above says that
$$\Shv_\Nilp(\Bun_G)=\Shv_{\on{HL}}(\Bun_G)$$
as subcategories of $\Shv(\Bun_G)$. 

\begin{rem}

In fact, for all but one statement in \cite{AGKRRV1}-\cite{AGKRRV3}, one needs the definition of 
$\Shv_\Nilp(\Bun_G)$ as $\Shv_{\on{HL}}(\Bun_G)$. 

\smallskip

The one place where one needs the definition via the singular support condition is the assertion that
the (a priori, fully faithful) functor
$$\Shv_\Nilp(\Bun_G)\otimes \Shv_\Nilp(\Bun_G)\to \Shv_{\Nilp\times \Nilp}(\Bun_G\times \Bun_G)$$
is an equivalence. 

\smallskip

The latter property exhibits a ``restricted" behavior of the category $\Shv_\Nilp(\Bun_G)$,
see \secref{sss:check global restr}. 

\end{rem} 

\ssec{The trace isomorphism}

In this subsection we explain how the subcategory $\Shv_\Nilp(\Bun_G)$ helps us achieve goal (a) in \secref{sss:goals}.

\sssec{}

Consider the embedding \eqref{e:Nilp as subcat}. It was conjectured in \cite{AGKRRV1} and subsequently proved in 
\cite{GR} that it admits a (colimit-preserving) right adjoint. Hence, by \secref{sss:funct Tr}, it induces a map
$$\Tr(\Frob_*,\Shv_\Nilp(\Bun_G))\to \Tr(\Frob_*,\Shv(\Bun_G)).$$

\sssec{}

The following result was established in \cite[Theorems 0.2.6 and 0.6.8]{AGKRRV3}:

\begin{thm} \label{t:Tr isom}
The composition
$$\Tr(\Frob_*,\Shv_\Nilp(\Bun_G))\to \Tr(\Frob_*,\Shv(\Bun_G)) \overset{\on{LT}^{\on{true}}}\longrightarrow \on{Funct}_c(\Bun_G(\BF_q))=:\on{Autom}(X,G)$$
is an isomorphism.
\end{thm} 

\begin{rem} 

\thmref{t:Tr isom} is powerful in that it establishes a direct connection between the geometric and classical theories: if we prove something interesting 
about the category $\Shv_\Nilp(\Bun_G)$, by taking the trace of Frobenius, we might be able to deduce interesting consequences for the space 
of automorphic functions. An example of this--the geometric Langlands equivalence--will be discussed in \secref{ss:restr geom Langlands}. 

\end{rem}

\ssec{The stack of local systems with restricted variation} \label{ss:LS restr}

In this subsection we let $k$ be an arbitrary ground field. 

\smallskip

As was mentioned in \secref{sss:fail goal}, we do not really know what the moduli stack\footnote{As long as we stay within \emph{algebraic}
(as opposed to analytic) geometry over $\ol\BQ_\ell$.} of $\cG$-local systems should be that 
could be relevant for the Langlands-type description of the entire category $\Shv(\Bun_G)$. However, an object that is suitable for
$\Shv_\Nilp(\Bun_G)$ was suggested in \cite{AGKRRV1}, under the name ``the stack of local systems with restricted variation",
denoted $\LS_\cG^{\on{restr}}(X)$. 

\sssec{}

The prestack $\LS_\cG^{\on{restr}}(X)$ is an object of derived algebraic geometry over the field of coefficients $\sfe$. Its values 
on a (connective) commutative $\sfe$-algebra $R$ is the ($\infty$)-groupoid of right t-exact symmetric monoidal functors
$$\Rep(\cG)\to R\mod(\Shv(X)).$$

It follows formally that any functor as above in fact takes values in the full subcategory
$$R\mod(\qLisse(X))\subset R\mod(\Shv(X)).$$

\sssec{}

The above definition works not just for \'etale $\ol\BQ_\ell$-adic sheaves, but for any constructible sheaf theory
(e.g., constructible Betti sheaves of $\sfe$-vector spaces when $k=\BC$ or holonomic D-modules when 
$k$ has characteristic $0$).

\smallskip

Let us compare this definition with the more familiar moduli stacks  
$$\LS_\cG^{\on{Betti}}(X) \text{ and } \LS_\cG^{\on{dR}}(X)$$
in the above two examples (see \cite[Sect. 4]{AGKRRV1} for more details). 

\sssec{}

The stack $\LS_\cG^{\on{Betti}}(X)$ attaches to $R$ the ($\infty$)-groupoid of right t-exact symmetric monoidal functors
$$\Rep(\cG)\to R\mod(\Shv^{\on{Betti}}_{\on{lisse}}(X)),$$
where $\Shv^{\on{Betti}}(X)$ is the category of all (i.e., not necessarily ind-constructible) Betti sheaves on $X$. 

\smallskip

It follows formally that any functor as above in fact takes values in the full subcategory
$$R\mod(\Shv^{\on{Betti}}_{\on{lisse}}(X)) \subset R\mod(\Shv^{\on{Betti}}(X)),$$
where $\Shv^{\on{Betti}}_{\on{lisse}}(X)$ is the category of all (i.e., not necessarily ind-constructible) lisse Betti sheaves on $X$. 

\smallskip

The stack $\LS_\cG^{\on{dR}}(X)$ attaches to $R$ the ($\infty$)-groupoid of right t-exact symmetric monoidal functors
$$\Rep(\cG)\to R\mod(\Dmod(X)).$$

\smallskip

In both cases we have the natural embeddings
$$\LS_\cG^{\on{restr,Betti}}(X)\subset \LS_\cG^{\on{Betti}}(X) \text{ and } \LS_\cG^{\on{restr,dR}}(X)\subset \LS_\cG^{\on{dR}}(X).$$ 

Moreover, one can describe their essential images as follows: $\LS_\cG^{\on{restr,Betti}}(X)$ (resp., $\LS_\cG^{\on{restr,dR}}(X)$)
is the disjoint union of formal completions of closed loci of $\LS_\cG^{\on{Betti}}(X)$ (resp., $\LS_\cG^{\on{dR}}(X) $) corresponding
to $\cG$-local systems with a given semi-simplification. 

\begin{rem}

In the case of Betti sheaves, one can give an even more explicit description of $\LS_\cG^{\on{restr,Betti}}(X)$:

\smallskip

Let $\LS_\cG^{\on{coarse,Betti}}(X)$ (a.k.a. character variety) be the coarse moduli space of $\cG$-local systems on $X$, i.e.,
$$\LS_\cG^{\on{coarse,Betti}}(X)=\Spec(\Gamma(\LS_\cG^{\on{Betti}}(X),\CO_{\LS_\cG^{\on{Betti}}(X)})).$$

We have a tautological map
$$\LS_\cG^{\on{Betti}}(X)\to \LS_\cG^{\on{coarse,Betti}}(X).$$

For an affine scheme $Z$ over $\sfe$, let $Z^{\text{take\,apart}}$ be the formal scheme equal to the disjoint union 
of formal completions of closed points of $Z$. Then we have a pullback square
$$
\CD
\LS_\cG^{\on{restr,Betti}}(X) @>>> \LS_\cG^{\on{Betti}}(X) \\
@VVV @VVV \\
(\LS_\cG^{\on{coarse,Betti}}(X))^{\text{take\,apart}} @>>> \LS_\cG^{\on{coarse,Betti}}(X).
\endCD
$$

\end{rem} 

\sssec{}

In \cite{AGKRRV1}, some basic properties of $\LS_\cG^{\on{restr}}(X)$ have been established. The most relevant 
for us is that it is a \emph{formal algebraic stack} locally almost of finite type (see {\it loc. cit.}, Theorem 1.4.5).
Moreover, it is \emph{quasi-smooth}\footnote{A.k.a., derived
locally complete intersection.}.  

\smallskip

In particular, the theory of singular support of coherent sheaves, developed in \cite{AG}, is applicable to $\LS_\cG^{\on{restr}}(X)$
(see \cite[Sect. 21.1]{AGKRRV1}). Thus, we can consider the full category
$$\IndCoh_\Nilp(\LS_\cG^{\on{restr}}(X))\subset \IndCoh(\LS_\cG^{\on{restr}}(X)),$$
see \cite[Sect. 11.1.5]{AG}. 

\begin{rem}
The category $\IndCoh_\Nilp(\LS_\cG^{\on{restr}}(X))$ should be thought of as a variant of the usual 
category $\QCoh(\LS_\cG^{\on{restr}}(X))$ that has to do with the fact that the formal stack $\LS_\cG^{\on{restr}}(X)$ is only quasi-smooth,
and not smooth. 
\end{rem}

\ssec{Restricted (global, unramified) geometric Langlands} \label{ss:restr geom Langlands} 

\sssec{}

The following result, which explains the relevance of $\LS_\cG^{\on{restr}}(X)$ to $\Shv_\Nilp(\Bun_G)$,
was established in \cite[Theorem 14.3.2]{AGKRRV1}:

\begin{thm} \label{t:sped decomp global}
There exists an action of $\QCoh(\LS_\cG^{\on{restr}}(X))$ on $\Shv_\Nilp(\Bun_G)$ so that the functor
$\on{H}_V$ of \eqref{e:Hecke functor} acts as tensor product by 
$$\CE_V\in \QCoh(\LS_\cG^{\on{restr}}(X))\otimes \qLisse(X),$$
the latter being the tautological object associated to $V\in \Rep(\cG)$.
\end{thm} 

\sssec{}

We are now ready to state the restricted version of the geometric Langlands conjecture:

\begin{conj} \label{c:restr GLC}
There exists an equivalence
$$\Shv_\Nilp(\Bun_G) \overset{\BL_G^{\on{restr}}}\simeq \IndCoh_\Nilp(\LS_\cG^{\on{restr}}(X)),$$
compatible with the action of $\QCoh(\LS_\cG^{\on{restr}}(X))$ on the two sides.
\end{conj}

The above \conjref{c:restr GLC} is close to be a theorem, but it is not one yet. In the paper \cite{GR}, 
\conjref{c:restr GLC} was established when $k$ has characteristic $0$, and when $G=GL_n$ over any ground field. 

\sssec{}

In fact, in \cite{GR} the following partial result in the direction of \conjref{c:restr GLC} has been established:

\begin{thm} \label{t:restr GLC}
There exists an equivalence
$$\Shv_\Nilp(\Bun_G) \overset{\BL_G^{\on{restr}}}\simeq \IndCoh_\Nilp({}'\!\LS_\cG^{\on{restr}}(X)),$$
compatible with the action of $\QCoh(\LS_\cG^{\on{restr}}(X))$ on the two sides, where 
$'\!\LS_\cG^{\on{restr}}(X)$ is the disjoint union of some of the connected components of 
$\LS_\cG^{\on{restr}}(X)$.
\end{thm}

\sssec{}

Thus \conjref{c:restr GLC} (resp., \thmref{t:restr GLC}) achieves goal (2) from \secref{sss:goals}.

\ssec{Consequences for automorphic functions} \label{ss:Tr of GLC}

Now with goals (1) and (2) achieved, we can relate the space $\on{Autom}(X,G)$ to the automorphic side.

\smallskip

In what follows we will assume the validity of \conjref{c:restr GLC}. The actual proven statements could
be obtained from the ones below by substituting $\LS_\cG^{\on{restr}}(X)\rightsquigarrow {}'\!\LS_\cG^{\on{restr}}(X)$.

\sssec{}

We once again assume that $k=\ol\BF_q$, but $G$ and $X$ are defined over $\BF_q$. The geometric Frobenius
on $X$ defines an automorphism, also denoted $\Frob$, of the prestack $\LS_\cG^{\on{restr}}(X)$. Consider
the resulting endofunctor (in fact, a self-equivalence) $\Frob_*$ of $\IndCoh_\Nilp(\LS_\cG^{\on{restr}}(X))$.

\smallskip

Since the construction
of the functor $\BL_G^{\on{restr}}$ is canonical\footnote{It actually depends on the choice of a square root of the canonical
bundle on $X$. Thus, we assume having made this choice; moreover we assume that this choice is Frobenius-equivariant.},
it intertwines the $\Frob_*$ endofunctors on both sides. 

\smallskip

In particular, from \conjref{c:restr GLC}, we obtain an isomorphism of vector spaces
\begin{equation} \label{e:Tr isomorphism}
\Tr(\Frob_*,\Shv_\Nilp(\Bun_G)) \simeq \Tr(\Frob_*, \IndCoh_\Nilp(\LS_\cG^{\on{restr}}(X))).
\end{equation} 

\smallskip

Let us now calculate both sides in \eqref{e:Tr isomorphism}.

\sssec{}

Let $\CY$ be a prestack\footnote{Assumed to have a schematic diagonal.}
such that the category $\QCoh(\CY)$ is dualizable. Let $\phi$ be an endomorphism of 
$\CY$. Then an elementary calculation (see \secref{sss:Tr QCoh}) shows that there is a canonical isomorphism
$$\Tr(\phi_*,\QCoh(\CY))\simeq \Gamma(\CY^\phi,\CO_{\CY^\phi}),$$
where $\CY^\phi$ is the prestack of $\phi$-fixed points on $\CY$.

\smallskip

Let now $\CY$ be a prestack locally almost of finite type; in this case, the category $\IndCoh(\CY)$ is well-defined. 
Assume that $\IndCoh(\CY)$ is dualizable. Then a parallel computation shows that there is a canonical 
isomorphism
\begin{equation} \label{e:Tr IndCoh}
\Tr(\phi_*,\IndCoh(\CY))\simeq \Gamma(\CY^\phi,\omega_{\CY^\phi}),
\end{equation} 
where $\omega_\CZ$ is the dualizing sheaf on $\CZ$. 

\sssec{}

Denote
$$\LS_\cG^{\on{arithm}}(X):=(\LS_\cG^{\on{restr}}(X))^{\Frob}.$$

We stipulate that $\LS_\cG^{\on{arithm}}(X)$ to be the stack of arithmetic $\cG$-local systems on $X$.
By definition its $\ol\BQ_\ell$-points are $\cG$-local systems on $X$, equipped with a Weil structure.

\smallskip

It is shown in \cite[Sect. Theorem 24.1.4]{AGKRRV1} that the prestack $\LS_\cG^{\on{arithm}}(X)$ is actually a 
quasi-compact algebraic stack.

\begin{rem}

Note that whereas $\LS_\cG^{\on{restr}}(X)$ is quasi-smooth, the stack $\LS_\cG^{\on{arithm}}(X)$ is quasi-quasi-smooth\footnote{The terminology was
suggested by D.~Beraldo.}. This means that the cotangent fibers of $\LS_\cG^{\on{restr}}(X)$ have cohomological amplitude $[-2,1]$
(the positive direction is due to stackiness, and the negative is due to the singularities). 

\smallskip

In particular, it is \emph{not} eventually coconnective, i.e., the structure sheaf has cohomologies in infinitely many 
negative cohomological degrees.

\end{rem} 

\sssec{}

By \eqref{e:Tr IndCoh}, we have:
\begin{equation} \label{e:Tr IndCoh LS}
\Tr(\Frob_*,\IndCoh(\LS_\cG^{\on{restr}}(X))\simeq \Gamma(\LS_\cG^{\on{arithm}}(X),\omega_{\LS_\cG^{\on{arithm}}(X)}).
\end{equation} 

In addition, we have:

\begin{prop}[{GR}, Theorem 1.5.2, quoting \cite{BLR}] \label{p:Nilp doesnt matter}
The inclusion 
$$\IndCoh_\Nilp(\LS_\cG^{\on{restr}}(X))\hookrightarrow \IndCoh(\LS_\cG^{\on{restr}}(X))$$
induces an isomorphism on $\Tr(\Frob_*,-)$.
\end{prop}

Combining \eqref{e:Tr IndCoh LS} with \propref{p:Nilp doesnt matter}, we obtain an isomorphism
\begin{equation} \label{e:Tr IndCoh LS Nilp}
\Tr(\Frob_*,\IndCoh_\Nilp(\LS_\cG^{\on{restr}}(X))\simeq \Gamma(\LS_\cG^{\on{arithm}}(X),\omega_{\LS_\cG^{\on{arithm}}(X)}).
\end{equation} 

\sssec{}

We now combine \eqref{e:Tr IndCoh LS Nilp}, \eqref{e:Tr isomorphism} and \thmref{t:Tr isom}. We obtain:

\begin{cor} \label{c:Autom descr}
There exists a canonical isomorphism
$$\on{Autom}(X,G)\simeq \Gamma(\LS_\cG^{\on{arithm}}(X),\omega_{\LS_\cG^{\on{arithm}}(X)}).$$
\end{cor} 

The isomorphism of \corref{c:Autom descr} can be viewed as an answer to Langlands' original question as
to the relationship between automorphic functions and $\cG$-local systems (in the unramified case).

\smallskip

In his talk, Sam Raskin will explain how to use \corref{c:Autom descr} to extract some desired information
about the space of automorphic functions, e.g., the Ramanujan conjecture\footnote{This is work-in-progress of
S.~Raskin, V.~Lafforgue and the author.} (in the unramified situation). 

\begin{rem} 

Note that, being the stack of fixed points, $\LS_\cG^{\on{arithm}}(X)$ is canonically Calabi-Yau. This
implies that over the quasi-smooth locus of $\LS_\cG^{\on{arithm}}(X)$, the dualizing sheaf is canonically 
isomorphic to the structure sheaf. 

\smallskip

However, the two are wildly different over the non quasi-smooth locus. In fact, the difference between
the two can be seen as the source of the non-temperedness phenomena. This will also be addressed
in Sam Raskin's talk.

\end{rem} 

\section{The local theory}  \label{s:local}

The contents of this section constitute joint work-in-progress, joint with S.~Raskin, and some parts are joint
also with other groups of mathematicians, to be named individually in the relevant subsections. 

\smallskip

We will propose a set of conjectures about the \emph{local} Langlands theory. There will be three such main conjectures:

\begin{itemize}

\item A 2-categorical local equivalence, \conjref{c:local GLC} (loosely, a counterpart of \conjref{c:restr GLC});

\item A local trace conjecture, \conjref{c:local Tr conj} (loosely, a counterpart of \thmref{t:Tr isom});

\item A 1-categorical local equivalence, \conjref{c:1-cat local Langlands} (loosely, a counterpart of \corref{c:Autom descr}).

\end{itemize}

The last conjecture is a formal consequence of the previous two, has been first proposed by X.~Zhu as \cite[Conjecture 4.6.4]{Zhu1}. 
and is closely related to the conjecture of L.~Fargues and P.~Scholze, \cite[Conjecture I.10.2]{FS}. 

\smallskip

We should also remark that the conjectures stated in this section seem quite tractable (and may already be theorems), mostly
thanks to a recent progress achieved in \cite{DVY}. 

\ssec{Categorical representations of groups--first attempt}

\sssec{}

Note that in the passage ``global" $\Rightarrow$ ``local" in the classical theory, one raises the level
of categoricity by one:

\smallskip

Globally we study the \emph{vector space} of automorphic functions (with the various operators acting on it),
whereas locally we study the \emph{category} of representations of a given $p$-adic group (or an enhanced
version of that--the category introduced by Fargues-Scholze in \cite{FS}). 

\smallskip

Note also that the passage ``classical" $\Rightarrow$ ``geometric" also raises the level of categoricity by one:

\smallskip

In the global theory, classically we study the \emph{vector space} of automorphic functions, whereas geometrically 
we study the \emph{category} of sheaves on $\Bun_G$ (or variants thereof). 

\smallskip

Thus, both of the above patterns show that the object of study the local geometric theory should be a \emph{2-category}.
We stipulate that the 2-category in question should be that of \emph{categorical representations of the loop group} $\fL(G):=G\ppart$.
Over the next few subsections we will explain what we mean by this. 

\smallskip

We start with the case when, instead of the loop group, we are dealing with an algebraic group $H$ of finite type.

\sssec{}

Let us first place ourselves in the context of D-modules. In this case, we consider the category $\Dmod(H)$, which has
a natural structure of \emph{commutative Hopf category}, i.e., this is a commutative Hopf algebra object in the 
symmetric monoidal 2-category $\DGCat$ (see \cite[Chapter 1, Sect. 10]{GaRo}):

\smallskip

The commutative algebra structure is given by \emph{pointwise tensor product}, and the coalgebra structure is given by pullback
along the multiplication map
$$\Dmod(H)\overset{\on{mult}^!}\to \Dmod(H\times H).$$

Here we are using a crucial fact that for a scheme $Y$ (and in fact for any reasonable prestack), the functor
\begin{equation} \label{e:cat Kunneth}
\Dmod(Y)\otimes \Dmod(Y')\overset{\boxtimes}\to \Dmod(Y\times Y')
\end{equation}
given by \emph{external tensor product}, is an equivalence for any prestack $Y'$. 
In what follows, we will refer to \eqref{e:cat Kunneth} 
as the \emph{categorical K\"unneth formula}. 

\smallskip

Given the above structure on $\Dmod(H)$, we define $H\on{-Cat}$ to be the 2-category 
$$\Dmod(H)\comod(\DGCat).$$

\smallskip

The commutative Hopf algebra structure on $\Dmod(H)$ defines on $H\on{-Cat}$ a symmetric monoidal
structure. 

\begin{rem}  \label{r:D(H) dual}

As is the case for any scheme of finite type, Verdier duality makes 
$\Dmod(H)$ self-dual as an object of $\DGCat$. Hence, we can regard $\Dmod(H)$
also as a \emph{cocommutative Hopf algebra} in $\DGCat$. 

\smallskip

The 2-category $H\on{-Cat}$ can be equivalently thought as \emph{modules} for 
$\Dmod(H)$ with respect to the above algebra structure. 

\end{rem} 

\sssec{}

We now return to the context of $\ol\BQ_\ell$-adic sheaves on schemes over $k$,
and we wish to define what $H\on{-Cat}$ is.

\smallskip

Naively, we might attempt what we just did in the case of D-modules. We start with the
category $\Shv(H)$, and... 

\smallskip

The problem is that the analog of the categorical K\"unneth formula \eqref{e:cat Kunneth}
fails for $\Shv(-)$. Namely, we have a fully faithful functor
\begin{equation} \label{e:cat Kunneth Shv}
\Shv(Y)\otimes \Shv(Y')\overset{\boxtimes}\to \Shv(Y\times Y'),
\end{equation}
but it is in general not an equivalence. In fact \eqref{e:cat Kunneth Shv} is \emph{never} an equivalence if
both $Y$ and $Y'$ are positive-dimensional schemes of finite type.

\sssec{} \label{sss:restr stack} 

For future reference, we give the following definition: a prestack $\CY$ is said to be \emph{restricted} if 
the functor
$$\Shv(\CY)\otimes \Shv(Y')\overset{\boxtimes}\to \Shv(\CY\times Y')$$
is an equivalence for any $Y\in \Sch^{\on{aff}}$.

\smallskip

An algebraic stack that has finitely many isomorphism classes of geometric points is restricted. An example of
such is 
$$\CY=N\backslash G/B.$$

\sssec{}

The failure of \eqref{e:cat Kunneth Shv} to be an equivalence means that we cannot define a structure of
coalgebra on $\Shv(H)$ using pullback along $H\times H\overset{\on{mult}}\to H$ ``as-is". 

\smallskip

One can, however, try the dual approach as in Remark \ref{r:D(H) dual}: i.e., we can treat $\Shv(H)$ as
a \emph{monoidal category} with respect to convolution. In this way we obtain \emph{a} definition of
categorical representations; we denote the resulting 2-category by $H\on{-Cat}^{\on{naive}}$.

\smallskip

However, we have not really resolved the issue. For example, $H\on{-Cat}^{\on{naive}}$ defined in the above 
way does \emph{not} have a symmetric monoidal structure compatible with the forgetful functor to $\DGCat$. 

\sssec{}

To produce a better-behaved object, we will have to change the paradigm: instead of looking for actions on $H$
on DG categories, we will replace the latter by a different symmetric monoidal 2-category, denoted $\AGCat$,
whose definition was suggested by V.~Drinfeld some 20 years. 

\ssec{The $(\infty,2)$-category $\AGCat$}

The contents of this and the next subsections constitute a joint project with N.~Rozenblyum and Y.~Varshavsky, \cite{GRV}.

\smallskip

The complete definitions, formulated in the proper homotopy-theoretic language, will be supplied in {\it loc. cit.} 
Here we give a brief digest.   

\sssec{}

An object $\ul\bC\in \AGCat$ is an assignment 
$$S\in \Sch^{\on{aff}} \mapsto \ul\bC_S\in \Shv(S)\mod(\DGCat),$$
where we regard $\Shv(S)$ as a monoidal category with respect to the !-tensor product.

\smallskip

The above assignment should be functorial in the following sense:

\begin{itemize}

\item For a map $f:S_1\to S_2$, we should have $\Shv(S_2)$-linear functors
$$f^!:\ul\bC_{S_2}\to \ul\bC_{S_1} \text{ and } f_*:\ul\bC_{S_1}\to \ul\bC_{S_2};$$

\item The above functors should be compatible with compositions of morphisms 
in the obvious sense;

\item For a Cartesian diagram of affine schemes
$$
\CD
S'_1 @>{g_1}>> S_1 \\
@V{f'}VV @VV{f}V \\
S'_2 @>{g_2}>> S_2,
\endCD
$$
we should be given an isomorphism
$$g_2^!\circ f_*\simeq f'_*\circ g_1^!$$
as functors $\ul\bC_{S_1}\to \ul\bC_{S'_2}$.

\end{itemize}

The above data must satisfy higher compatibilities, which can be formulated using the formalism
of \cite[Chapter 7]{GaRo}.

\smallskip

Morphisms in $\AGCat$ are defined naturally: for a pair of objects $\ul\bC'$ and $\ul\bC''$, we define the 
DG category of 1-morphisms between them, to have as objects compatible collections of $\Shv(S)$-linear functors 
$$\ul\bC'_S\to \ul\bC''_S, \quad S\in \Sch^{\on{aff}}$$
compatible with the operations $f^!$ and $f_*$, etc. 

\sssec{}

To a prestack $\CY$ one attaches an object of $\AGCat$, denoted $\ul\Shv(\CY)$, defined by
$$\ul\Shv(\CY)_S:=\Shv(\CY\times S),$$
with the operations $f^!$ and $f_*$ defined naturally. 

\sssec{}

We have a naturally defined fully faithful functor
$$\on{emb.restr}:\DGCat\to \AGCat$$
that assigns to $\bC\in \DGCat$ the object
$$S\mapsto \bC\otimes \Shv(S)$$
with the operations $f^!$ and $f_*$ defined naturally. 

\smallskip

The functor $\on{emb.restr}$ admits a right adjoint, given by
$$\ul\bC\mapsto \ul\bC_{\on{pt}}.$$

\sssec{}

In what follows we will say that an object of $\AGCat$ is \emph{restricted} if it lies in the essential image of 
the functor $\on{emb.restr}$. 

\smallskip

For example, a prestack $\CY$ is restricted in the sense of \secref{sss:restr stack} if and only if 
$\ul\Shv(\CY)$ is restricted as an object of $\AGCat$. 

\sssec{}

A key construction that makes $\AGCat$ do the job for us is that of the symmetric monoidal
structure on it. One shows that $\AGCat$ possesses a unique symmetric monoidal structure, such that:

\begin{itemize}

\item It is compatible with colimits;

\item The functor $\on{emb.restr}$ is symmetric monoidal;

\item For $S_1,S_2\in \Sch^{\on{aff}}$, we have 
\begin{equation} \label{e:cat Kunneth AGCat aff}
\ul\Shv(S_1)\otimes \ul\Shv(S_2)=\ul\Shv(S_1\times S_2), \quad S_1,S_2\in \Sch^{\on{aff}}.
\end{equation} 

\end{itemize} 

The unit of this symmetric monoidal structure is the object
$$\ul\Vect:=\ul\Shv(\on{pt})\simeq \on{emb.restr}(\Vect).$$

\sssec{}

One shows that for a reasonably large class of prestacks, the naturally defined map
\begin{equation} \label{e:cat Kunneth AGCat stacks}
\ul\Shv(\CY_1)\otimes \ul\Shv(\CY_2)\to \ul\Shv(\CY_1\times \CY_2)
\end{equation} 
is still an isomorphism.

\sssec{}

Note that since the functor $\on{emb.restr}$ is symmetric monoidal, its right adjoint possesses a 
\emph{right-lax} symmetric monoidal structure. Evaluating on the objects 
$\ul\Shv(\CY_1),\ul\Shv(\CY_2)\in \AGCat$, we obtain a 1-morphism in $\DGCat$ (i.e., a functor)
\begin{multline} \label{e:cat Kunneth DGCat stacks again}
\Shv(\CY_1)\otimes \Shv(\CY_2)\simeq
\on{emb.restr}^R(\ul\Shv(\CY_1))\otimes \on{emb.restr}^R(\ul\Shv(\CY_2)) \to \\
\to \on{emb.restr}^R(\ul\Shv(\CY_1)\otimes \ul\Shv(\CY_2)) \to \on{emb.restr}^R(\ul\Shv(\CY_1\times \CY_2))=\Shv(\CY_1\times \CY_2),
\end{multline}
which is the usual $\boxtimes$ functor. 

\smallskip

Note that \eqref{e:cat Kunneth DGCat stacks again} is typically not an isomorphism, even when \eqref{e:cat Kunneth AGCat stacks} is. 

\sssec{} \label{sss:trace in AGCat} 

Let $\CY$ be a quasi-compact algebraic stack. Under some mild assumptions 
on $\CY$, the map \eqref{e:cat Kunneth AGCat stacks} is an isomorphism for $\CY_1=\CY_2=\CY$. This formally
implies that $\ul\Shv(\CY)$ is dualizable (and self-dual) as an object of $\AGCat$. 

\smallskip

Let $\phi$ be an endomorphism of $\CY$ (assumed to be schematic, so that the functor $\phi_*$ is well-defined).
In this case, it makes sense to consider
$$\Tr_{\AGCat}(\phi_*,\ul\Shv(\CY))\in \End(\ul\Vect)\simeq \Vect.$$

A formal computation shows:
\begin{equation} \label{e:trace in AGCat} 
\Tr_{\AGCat}(\phi_*,\ul\Shv(\CY)) \simeq \on{C}^\cdot_\blacktriangle(\CY^\phi,\omega_{\CY^\phi}),
\end{equation} 
where $\on{C}^\cdot_\blacktriangle(-)$ is the functor of \emph{renormalized cochains}, see \cite[Sect. A.2.4]{AGKRRV2}.

\smallskip

In particular, the category $\AGCat$ helps us turn the map \eqref{e:Tr Frob on sheaves qc} into an isomorphism:
for $k=\ol\BF_q$ with $\CY$ defined over $\BF_q$ and $\phi=\Frob$, the map \eqref{e:Tr Frob on sheaves qc} 
factors as
$$\Tr_{\DGCat}(\Frob_*,\Shv(\CY)) \overset{\on{emb.restr}}\longrightarrow \Tr_{\AGCat}(\Frob_*,\ul\Shv(\CY))\simeq
\on{C}^\cdot_\blacktriangle(\CY^{\Frob},\omega_{\CY^{\Frob}})\simeq \on{Funct}(\CY(\BF)_q).$$

\smallskip

The failure of the initial map \eqref{e:Tr Frob on sheaves qc} to be an isomorphism stems from the failure of 
\eqref{e:cat Kunneth DGCat stacks again} to be an isomorphism. 

\ssec{Categorical representations of groups--the definition}

Having introduced the symmetric monoidal 2-category, we can now define what we mean by a categorical $H$-representation.

\sssec{} \label{sss:H-cat}

By construction, we can view $\ul\Shv(H)$ as a commutative Hopf algebra in $\AGCat$. We define the 2-category
$H\on{-Cat}$ to be 
$$\ul\Shv(H)\comod(\AGCat).$$

\sssec{}

Thus, explicitly, an object of $H\on{-Cat}$ is an object 
$$\ul\bC=\{S\mapsto \ul\bC_S\}$$
of $\AGCat$, equipped with the datum of functors
$$\on{co-act}_S:\ul\bC_S\to \ul\bC_{S\times H}$$
that are associative in the sense that the diagrams 
$$
\CD
\ul\bC_S @>{\on{co-act}_S}>> \ul\bC_{S\times H} \\ 
@V{\on{co-act}_S}VV @VV{\on{co-act}_{H\times S}}V \\ 
\ul\bC_{S\times H} @>{(\on{id}_S\times \on{mult}_H)^!}>> \ul\bC_{S\times H\times H}
\endCD
$$
commute. 

\sssec{Example}

Let $\CY$ be a prestack equipped with an action of $H$. Then $\ul\Shv(Y)$ naturally upgrades to an object
of $H\Cat$. 

\sssec{} \label{sss:properties H-Cat}

The commutative algebra structure on $\ul\Shv(H)$ defines on $H\on{-Cat}$ a symmetric monoidal structure, compatible
with the forgetful functor 
\begin{equation} \label{e:oblv H}
\oblv_H:H\on{-Cat}\to \AGCat.
\end{equation} 

\smallskip

An object of $\ul\Shv(H)$ is dualizable if and only if its image under $\oblv_H$
is dualizable, and $\oblv_H$ commutes with duality.

\smallskip

A somewhat non-trivial fact is that the functor
$$\on{inv}_H:H\on{-Cat}\to \AGCat, \quad \ul\bC\mapsto \ul\bC^H$$
also commutes with duality (see \cite[Sect. B]{Ga(W)}). 

\smallskip

By construction, $H\on{-Cat}$ is tensored over $\AGCat$. One shows that it is \emph{rigid} as a symmetric
monoidal 2-category over $\AGCat$. In particular, it is dualizable and self-dual as a module over $\AGCat$.

\sssec{} \label{sss:naively restricted}

One can ask the following question (if one still desperately wants to consider actions of $H$ on DG categories, rather than 
on objects of the weird $\AGCat$): 

\smallskip

Can one describe the subcategory 
\begin{equation} \label{e:naively restricted}
H\on{-Cat}\underset{\AGCat}\times \DGCat,
\end{equation} 

where:

\begin{itemize}

\item $H\on{-Cat}\to \AGCat$ is the forgetful functor;

\item $\DGCat\to \AGCat$ is the embedding $\on{emb.restr}$?

\end{itemize}

Unwinding, these are DG categories $\bC$, equipped with a coaction functor 
$$\on{co-act}:\bC\to \bC\otimes \Shv(H),$$
such that the following diagram commutes:
$$
\CD
\bC @>{\on{co-act}}>> \bC\otimes \Shv(H) \\ 
& & @VV{\on{co-act}\otimes \on{Id}_{\Shv(H)}}V \\
@V{\on{co-act}}VV   \bC\otimes \Shv(H)\otimes \Shv(H) \\ 
& & @VV{\on{Id}_\bC \otimes (\boxtimes)}V \\
\bC\otimes \Shv(H) @>{\on{Id}_\bC\otimes \on{mult}_H^!}>> \bC\otimes \Shv(H\times H). 
\endCD
$$

One can classify all such $\bC$, but as we shall presently see, there will be relatively few such categories.

\sssec{} \label{sss:naively restricted expl}

Let $\chi\in \Shv(H)$ be a character sheaf on $H$, i.e.,
$$\on{mult}^!(\chi)\simeq \chi\boxtimes \chi.$$

Let $\Shv(H)_\chi\subset \Shv(H)$ be the full subcategory consisting of objects all of whose cohomologies
are extensions of copies of $\chi$. The operation $\on{mult}^!$ defines on $\Shv(H)_\chi$ a structure of comonoidal category.

\smallskip

It is clear that the functor $\on{emb.restr}$ extends to a (fully faithful) functor
$$\on{emb.restr}_\chi:\Shv(H)_\chi\comod(\DGCat)\to H\on{-Cat}.$$
%and further
%\begin{equation} \label{e:chi part}
%\Shv(H)_\chi\comod(\AGCat)\simeq \Shv(H)_\chi\comod(\DGCat)\underset{\DGCat}\otimes \AGCat
%\to H\on{-Cat}
%\end{equation} 

\smallskip

We have:

\begin{prop} \label{p:naively restricted}
Suppose that $H$ is connected. Then 
any object in \eqref{e:naively restricted} splits (canonically) as a direct sum of objects, each of
which lies in the essential image of the functor $\on{emb.restr}_\chi$ for some $\chi$.
\end{prop} 

\sssec{}

We now consider the case of the loop group. Since $\fL(G)$ is of infinite type, one needs to specify what one means 
by the category $\Shv(\fL(G))$. This can be done by the recipe of \cite[Sect. C.3]{Ga(W)}. In particular, we obtain a well-defined 
object 
$$\ul\Shv(\fL(G))\in \AGCat,$$
and the 1-morophism
$$\ul\Shv(\fL(G))\otimes \ul\Shv(\fL(G))\to \ul\Shv(\fL(G)\times \fL(G))$$
is an isomorphism.

\smallskip

Therefore, we can follow the recipe of \secref{sss:H-cat} and define the 
2-category
$$\fL(G)\Cat.$$

The properties mentioned in \secref{sss:properties H-Cat} apply in the present situation as well (see \cite[Sect. D]{Ga(W)}). 

\sssec{Examples}

The objects $\ul\Shv(\Gr_G)$, $\ul\Shv(\on{Fl}^{\on{aff}}_G)$ and $\ul\Shv(\fL(G))$ itself 
all naturally upgrade to objects of $\fL(G)\Cat$.

\ssec{Restricted representations}

The notion of \emph{restricted categorical representation} for a finite-dimensional reductive group 
has been developed in \cite{GRV}. In the case of loop groups, the quest for this notion has been 
a joint project with D.~Kazhdan, S.~Raskin and Y.~Varshavsky. However, a recent work of
G.~Dhillon, Y.~Varshavsky and D.~Yang (\cite{DVY}) seems to have settled all the basic questions
related to this notion (under the assumption that $p=\on{char}(k)>|W|$). 
 
\sssec{}

We now have the 2-category $\fL(G)\Cat$ up and running, but... this is still not good enough:

\smallskip

We have worked hard to define $\fL(G)\Cat$, and it came out to be a category tensored over $\AGCat$.
However, if our goal is to establish a Langlands-type equivalence between $\fL(G)$-categories and some 
2-category on the spectral side, the latter will only be tensored over $\DGCat$. 

\smallskip

In other words, $\fL(G)\Cat$ is too big in order be equivalent to something on the spectral side, and we need
to trim it down, i.e., find a full subcategory
\begin{equation} \label{e:restricted LG}
\fL(G)\Cat_{\on{restr}}\subset \fL(G)\Cat,
\end{equation}
which is tensored over $\DGCat$ and can potentially be equivalent to a 2-category on the spectral side.

\smallskip

One could regard the passage 
$$\fL(G)\Cat\rightsquigarrow \fL(G)\Cat_{\on{restr}}$$
as loosely analogous\footnote{This analogy is less loose if one interprets $\Shv_\Nilp(\Bun_G)$ as $\Shv_{\on{HL}}(\Bun_G)$
(see \secref{sss:Hecke lisse}) and if one also gives a similar interpretation to $\fL(G)\Cat_{\on{restr}}$, which conjecturally does exist.} 
to the passage
$$\Shv(\Bun_G)\rightsquigarrow \Shv_\Nilp(\Bun_G)$$
in \secref{ss:BunG Nilp}. 

\smallskip

Note that the naive guess of taking $\fL(G)\Cat_{\on{restr}}$ to be \eqref{e:naively restricted} cannot possibly 
work: as was explained in \secref{sss:naively restricted expl}, the resulting 2-category will be too small. 

\smallskip

Before we define \eqref{e:restricted LG} for the loop group, we consider the case of a finite-dimensional
reductive group $G$.

\sssec{}

We define the full subcategory
\begin{equation} \label{e:fin dim restr}
G\Cat_{\on{restr}}\subset G\Cat
\end{equation} 
as follows. 

\smallskip

Let $\ul\bC$ be an object of $G\Cat$. We say that it is \emph{restricted} if the object
$$\ul\bC^N:=\on{inv}_N(\ul\bC)\in \AGCat$$
is restricted. 

\smallskip

This definition is justified by the following:

\begin{lem}
The operation
$$\ul\bC\mapsto \ul\bC^N, \quad G\Cat\to \AGCat$$
is conservative.
\end{lem}

\sssec{Examples}

A typical example of an object of $G\Cat$ that belongs to $G\Cat_{\on{restr}}$ is $\ul\Shv(G/B)$. 
Note, however, that $\ul\Shv(G/B)$ is \emph{not} restricted as an object of $\AGCat$, i.e, when
we forget the $G$-action. 

\smallskip

A generalization of this example is
\begin{equation} \label{e:chi twisted B}
\ul\Shv(G/B)^\chi:=(\ul\Shv(G/N)\otimes \ul\Vect^\chi)^T,
\end{equation} 
where $\ul\Vect^\chi$ is a copy of $\ul\Vect$ as an object of $\AGCat$, on which $T$ acts
via the character sheaf $\chi$. 

\smallskip

An example of an object of $G\Cat$ that does \emph{not} belong to $G\Cat_{\on{restr}}$ is $\ul\Shv(G)$. 

\sssec{}

Note that for $\ul\bC\in G\Cat$, the object $\ul\bC^N\in \AGCat$ naturally upgrades to an object $T\Cat$.  
Hence, if $\ul\bC$ is restricted, by \secref{sss:naively restricted expl}, $\ul\bC^N$ splits as a direct sum,
indexed by character sheaves on $T$. 

\smallskip

From here it follows that the objects \eqref{e:chi twisted B} generate $G\Cat_{\on{restr}}$. In other words,
$G\Cat_{\on{restr}}$ is equivalent to the category of modules (in $\DGCat$) over the monoidal category
$$\underset{\chi',\chi''}\oplus\, \Shv(B\backslash G/B)^{\chi',\chi''}.$$

\sssec{}

We now consider the case of $\fL(G)$. Below we give a provisional definition. Another  
definition (but one that is hopefully equivalent to the one we give), but one that gives much 
more information about the behavior of the restrictedness condition has been 
recently obtained in a work of G.~Dhillon, Y.~Varshavsky and D.~Yang, \cite{DVY}. 

\smallskip

For $\ul\bC\in \fL(G)\Cat$ and a character sheaf $\psi$ of the group $\fL(N)$ consider
\begin{equation} \label{e:gen Whit cond}
\ul\bC^{\fL(N),\psi}:=(\ul\bC\otimes \ul\Vect^\psi)^{\fL(N)}\in \AGCat.
\end{equation}

We shall say that $\ul\bC$ is \emph{restricted as an object of} $\fL(G)\Cat$ if the objects
\eqref{e:gen Whit cond} are restricted as objects of $\AGCat$.

\sssec{Examples}

The objects $\ul\Shv(\Gr_G)$, $\ul\Shv(\on{Fl}^{\on{aff}}_G)$ are restricted 
\emph{as objects of} $\fL(G)\Cat$ (but not as objects of $\AGCat$). 

\smallskip

A typical example of an object of $\fL(G)\Cat$ that is \emph{not} restricted is
$\ul\Shv(\fL(G))$. 

\sssec{} \label{sss:c:restr properties}

Here are some properties of the inclusion \eqref{e:restricted LG} that are non-obvious from the 
definition we gave, but that follow from the description in \cite{DVY}: 

\begin{conj} \label{sss:restr properties}
The functor
$$\iota:\AGCat\underset{\DGCat}\otimes \fL(G)\Cat_{\on{restr}}\to \fL(G)\Cat,$$
induced by \eqref{e:restricted LG} is fully faithful and admits a right adjoint. Moreover:

\smallskip

\noindent{\em(a)} The counit of the above adjunction itself admits a right adjoint. 

\smallskip

\noindent{\em(b)} The self-duality on $\fL(G)\Cat$ over $\AGCat$ induces (automatically, uniquely)
a self-duality of $\fL(G)\Cat_{\on{restr}}$ over $\DGCat$ so that $\iota^\vee\simeq \iota^R$.

\end{conj}

In concrete terms, point (b) above says that a dualizable object of $\fL(G)\Cat$ belongs to $\fL(G)\Cat_{\on{restr}}$
if and only if its dual does, and if $\ul\bC_1,\ul\bC_2\in \fL(G)\Cat_{\on{restr}}$, then the object
\begin{equation} \label{e:inv restr}
(\ul\bC_1\otimes \ul\bC_2)^{\fL(G)}\in \AGCat
\end{equation} 
actually belongs to the essential image of $\DGCat$ under $\on{emb.restr}$. 

\smallskip

We note that properties parallel to ones in \propref{sss:restr properties} are easy to show for a finite-dimensional
reductive group $G$ and the embedding \eqref{e:fin dim restr}. 

\ssec{The local trace isomorphism}

Recall that in the global unramified theory, the operation $\Tr(\Frob,-)$ realized a passage between the geometric and classical
theories, while lowering the categorical level by one.

\smallskip 

In this subsection we explain that a similar procedure is applicable for the local theory, and produces from the 2-category
$\fL(G)\mod$

\sssec{}

Let $\fC$ be a (dualizable) 2-category tensored over $\DGCat$ (see \secref{sss:2 cats}), and let $\fF$ be an endofunctor of $\fC$. 
In this case, we can consider 
$$\Tr_{\DGCat}(\fF,\fC)\in \DGCat,$$
see \secref{sss:Tr 2 cats}. 

\smallskip

Similarly, if $\fC$ is a (dualizable) 2-category tensored over $\AGCat$ and $\fF$ its endofunctor, we can consider
 $$\Tr_{\AGCat}(\fF,\fC)\in \AGCat.$$
 
 \sssec{}

Let $H$ be a finite-dimensional algebraic group, and let $\phi$ be its endomorphism. Precomposition with 
$\phi$ defines an endomorphism of the 2-category $H\Cat$, which we denote by the same symbol $\phi$. 

\smallskip

Consider the corresponding object
$$\Tr_{\AGCat}(\phi,H\Cat)\in \AGCat.$$

A calculation parallel to \eqref{e:trace in AGCat} shows:
\begin{equation} \label{e:Tr H Cat}
\Tr_{\AGCat}(\phi,H\Cat)\simeq \ul\Shv(H/\on{Ad}_\phi(H)),
\end{equation}
where $\on{Ad}_\phi(H)$ refers to the action of $H$ on itself by $\phi$-twisted conjugation.

\sssec{Example} \label{sss:Lang}

Let $k=\ol\BF_q$, and let $H$ be connected and defined over $\BF_q$. Take $\phi=\Frob$.

\smallskip

Note now that Lang's theorem identifies
$$H/\on{Ad}_{\Frob}(H)\simeq \on{pt}/H(\BF_q).$$ 

This implies in particular that the object 
$$\ul\Shv(H/\on{Ad}_{\Frob}(H))\in \AGCat$$
is restricted and identifies with 
$$\on{emb.restr}(\Shv(\on{pt}/H(\BF_q)))\simeq \on{emb.restr}(\Rep(H(\BF_q))).$$

\sssec{} \label{sss:cl of cat rep}

Let $\ul\bC$ be an object of $H\Cat$, assumed dualizable\footnote{Recall that this is the same as asking that the underlying object of $\AGCat$
be dualizable. Note also that at the present 2-categorical level, this condition is satisfied in most natural examples.}.
Let $\ul\bC$ be equipped with a datum of weak compatibility with $\phi$, i.e., a morphism
$$\ul\bC\overset{\alpha}\to \phi(\ul\bC).$$

We can think of $\alpha$ as an endomorphism of $\ul\bC$ as an object of $\AGCat$ that intertwines the given action of $H$ on $\ul\bC$ and
the one precomposed with $\phi$. 

\smallskip

Then, by \secref{sss:class}, we can associate to the pair $(\ul\bC,\alpha)$ a 1-morphism
$$\on{cl}(\ul\bC,\alpha)\in \Hom(\one_{\AGCat},\Tr_{\AGCat}(\phi,H\Cat))\simeq \Tr_{\AGCat}(\phi,H\Cat)_{\on{pt}}\simeq 
\Shv(H/\on{Ad}_\phi(H)).$$

Unwinding, we obtain that !-fibers of $\on{cl}(\ul\bC,\alpha)$, viewed as a sheaf on $H/\on{Ad}_\phi(H)$, are given by
$$h\mapsto \Tr(\alpha\circ h,\ul\bC_{\on{pt}}),$$
where we view both $h\in H$ and $\alpha$ as endomorphisms of $\ul\bC_{\on{pt}}$. 

\sssec{} \label{sss:DL}

We continue to be in the setting of \secref{sss:cl of cat rep}, with $H=G$ be a connected reductive group over 
$\ol\BF_q$. Assume that $G$ is defined and split over $\BF_q$; let $B$ be a Borel subgroup, also
assumed defined over $\BF_q$. 

\smallskip

Take $\ul\bC=\ul\Shv(G/B)^\chi$, where $\chi$ is a character 
sheaf on $T$, such that
$$\Frob^*(\chi)=w\cdot \chi$$
for an element $w$ of the Weyl group. 

\smallskip

We have a (canonically defined) intertwining functor 
$$T_w:\ul\Shv(G/B)^{w\cdot \chi}\to \ul\Shv(G/B)^\chi,$$
which is in fact an isomorphism between objects of $G\Cat$.

\smallskip

We take $\alpha$ to be the composition
$$\ul\Shv(G/B)^\chi \overset{\Frob}\simeq \ul\Shv(G/B)^{w\cdot \chi}\overset{T_w}\to \ul\Shv(G/B)^\chi.$$

Unwinding, we obtain that the resulting object
$$\on{cl}(\ul\Shv(G/B)^\chi,\alpha)\in \Shv(G/\on{Ad}_{\Frob}(G))\simeq \Rep(\cG(\BF_q))$$
is the Deligne-Lusztig representation attached to the pair $(\chi,w)$. 

\sssec{} \label{sss:2-cat trace map fin dim}

Let $G$ be as above, and consider the object
$$\Tr_{\DGCat}(\Frob,G\Cat_{\on{restr}})\in \DGCat.$$

The embedding
$$G\Cat_{\on{restr}}\to G\Cat$$
is linear with respect to $\on{emb.restr}:\DGCat\to \AGCat$
and hence defines a map
\begin{equation} \label{e:2 Cat Tr fin dim}
\on{emb.restr}(\Tr_{\DGCat}(\Frob,G\Cat_{\on{restr}}))\to \Tr_{\AGCat}(\Frob,G\Cat)
\end{equation}

The following is equivalent to \cite[Theorem 6.1.1]{Ete}:

\begin{thm} \label{t:Tr isom fin dim}
The map \eqref{e:2 Cat Tr fin dim} is an isomorphism.
\end{thm} 

Thus, combined with the computation \eqref{e:Tr H Cat} and \secref{sss:Lang}, we obtain an isomorphism
\begin{equation} \label{e:Tr on restr fin dim}
\Tr_{\DGCat}(\Frob,G\Cat_{\on{restr}})\simeq \Rep(G(\BF_q))
\end{equation}
as objects in $\DGCat$.

\sssec{}

We now consider the case of the loop group: 

\smallskip

We \emph{define} the object
\begin{equation}  \label{e:Shvs on Isoc}
\ul\Shv(\on{Isoc}_G)\in \AGCat
\end{equation} 
to be 
$$\ul\Shv(\fL(G))^{\on{Ad}_{\Frob}(\fL(G))}.$$

\smallskip

Denote
$$\Shv(\on{Isoc}_G):=\ul\Shv(\on{Isoc}_G)_{\on{pt}}\in \DGCat.$$

\smallskip

The computation leading to \eqref{e:Tr H Cat} is valid in this case as well and yields
\begin{equation} \label{e:Tr LG AGCat}
\Tr_{\AGCat}(\Frob,\fL(G)\Cat)\simeq \ul\Shv(\on{Isoc}_G)
\end{equation}
as objects of $\AGCat$. 

\sssec{} \label{sss:Isoc geom}

One can actually consider the quotient\footnote{The idea to consider $\Shv(\on{Isoc}_G)$ as a ``correct" enlargement of 
$\Rep(\fL(G)(\BF_q))$, with the motivation provided by \eqref{e:Tr LG AGCat}, has probably occurred to several people 
at various points in time. In particular, it occurred to V.~Lafforgue and the author around 2013.}
$$\on{Isoc}_G:= \fL(G)/\on{Ad}_{\Frob}(\fL(G))$$
as a prestack over $\ol\BF_q$; a comprehensive study of this object 
has been recently undertaken by X.~Zhu in \cite{Zhu2}.  

\smallskip

Being a group ind-scheme, $\fL(G)$ does not satisfy Lang's theorem. However, it does contain
$\on{pt}/\fL(G)(\BF_q)$ as a closed subfunctor, so that we have a pair of adjoint functors
\begin{equation} \label{e:unit stratum}
\iota_!:\Rep(\fL(G)(\BF_q))\rightleftarrows \Shv(\on{Isoc}_G):\iota^!.
\end{equation} 

\smallskip

One of the key results about it is that the prestack $\on{Isoc}_G$ is \emph{restricted}, i.e., the  
object \eqref{e:Shvs on Isoc} is restricted (see \cite[Corollary 3.71]{Zhu2}), and so
\begin{equation} \label{e:Isoc is restr}
\ul\Shv(\on{Isoc}_G)\simeq \on{emb.restr}(\Shv(\on{Isoc}_G)).
\end{equation} 

\smallskip

For future use, we note that from the definition of $\Shv(\on{Isoc}_G)$ via \eqref{e:Shvs on Isoc}, it follows
that this category is dualizable and self-dual (see \cite[Proposition 3.82]{Zhu2}).

\sssec{}

As in \secref{sss:cl of cat rep}, given an object $\ul\bC\in \fL(G)\Cat$ (satisfying an appropriate compactness condition)
and a map
$$\alpha:\ul\bC\to \Frob(\ul\bC),$$
we can consider its class
$$\on{cl}(\ul\bC,\alpha)\in \Shv(\on{Isoc}_G).$$

For specific examples of $(\ul\bC,\alpha)$, one can view the object $\on{cl}(\ul\bC,\alpha)$ as a loop group
generalization of Deligne-Lusztig representations. 

\sssec{}

We now state a local 2-categorical counterpart of \thmref{t:Tr isom} (which is also a loop group
analog of \thmref{t:Tr isom fin dim}). 

\smallskip

We consider $\fL(G)\Cat_{\on{restr}}$ as a 2-category tensored over $\DGCat$, and equipped with an endomorphism given by Frobenius.

\smallskip

Consider the object
$$\Tr_{\DGCat}(\Frob,\fL(G)\Cat_{\on{restr}})\in \DGCat.$$

As in \secref{sss:2-cat trace map fin dim}, the functor \eqref{e:restricted LG} induces a 1-morphism

\begin{equation} \label{e:local trace map}
\on{emb.restr}(\Tr_{\DGCat}(\Frob,\fL(G)\Cat_{\on{restr}}))\to \Tr_{\AGCat}(\Frob,\fL(G)).
\end{equation}

We propose: 

\begin{conj} \label{c:local Tr conj}
The map \eqref{e:local trace map} is an isomorphism.
\end{conj}

Combining with \eqref{e:Tr LG AGCat} and \eqref{e:Isoc is restr}, \conjref{c:local Tr conj} implies
that the (a priori defined functor)
\begin{equation} \label{e:local trace conj}
\Tr_{\DGCat}(\Frob,\fL(G)\Cat_{\on{restr}})\to \Shv(\on{Isoc}_G)
\end{equation} 
is an isomorphism.

\sssec{}

Recently, in a joint work with G.~Dhillon, A.~Eteve, S.~Raskin, Y.~Varshavsky and D.~Yang a substantial progress
has been made towards \conjref{c:local Tr conj}, and independently by C.~Chan, T.~Kaletha and X.~Zhu by a different method.
(Un)surprisingly, this statement ended up being much simpler than its its global counterpart, \thmref{t:Tr isom}. 

\ssec{The restricted local geometric Langlands conjecture}

In this subsection we connect our 2-category of categorical representations of $\fL(G)$ to the spectral side. 

\sssec{}

We define the prestack (of $\cG$-local systems on the formal punctured disc $\ocD$), denoted
$\LS^{\on{restr}}_\cG(\ocD)$, by the same recipe as in \secref{ss:LS restr}, where instead of the category $\qLisse(X)$,
we consider
$$\qLisse(\ocD):=\on{Inert}\mod,$$
where $\on{Inert}$ is the inertia group associated with the field $k\ppart$. 

\smallskip

Like in the case of $\LS_\cG(X)$, the prestack $\LS^{\on{restr}}_\cG(\ocD)$ is a formal algebraic stack. Moreover,
since $\qLisse(\ocD)$ has cohomological dimension $1$, the stack $\LS^{\on{restr}}_\cG(\ocD)$ is actually formally smooth
(and not just quasi-smooth). 

\sssec{}

Naively, the \emph{restricted} local geometric conjecture would be an equivalence between 
$\fL(G)\Cat_{\on{restr}}$ and the 2-category of modules over $\QCoh(\LS^{\on{restr}}_\cG(\ocD))$ in $\DGCat$.
However, this naive guess must be undergo a correction 
$$\QCoh \rightsquigarrow \IndCoh,$$
which is a 2-categorical counterpart of the correction
$$\QCoh(\LS_\cG(X))\rightsquigarrow \IndCoh_\Nilp(\LS_\cG(X))$$
in the global unramified case. We shall presently explain the nature of this correction, following D.~Arinkin.\footnote{An alternative approach
was subsequently developed by G.~Stefanich in \cite{Ste}.} 

\sssec{} \label{sss:Arinkin}

Let $S$ be a smooth affine scheme. Denote
$$\tQCoh(S):=\QCoh(S)\mod(\DGCat).$$  

Let now $\CN \subset T^*(S)$ be a conical sub-Lagrangian subset. We will define a 2-category 
$\tIndCoh_\CN(S)$, which will agree (resp., contain) $\tQCoh(S)$ if $\CN$ equals (resp., contains)
the zero-section. 

\smallskip

Let $\wt{S}$ be a scheme equipped with a quasi-smooth and proper map 
$$f:\wt{S}\to S,$$
such that $\CN$ is contained in the locus
$$\{(s,\xi)\in T^*_s(S) \,|\, \exists \,\wt{s}\in \wt{S} \text{ such that } f(\wt{s})=s \text{ and } 
\xi\in \on{ker}(df:T_s^*(S)\to H^0(T_{\wt{s}}^*(\wt{S})))\}.$$

To this $\CN$ we will associate a monoidal localization of $\IndCoh(\wt{S}\underset{S}\times \wt{S})$, denoted
$\IndCoh_{\wt\CN}(\wt{S}\underset{S}\times \wt{S})$, and we will define 
$$\tIndCoh_\CN(S):=\IndCoh_{\wt\CN}(\wt{S}\underset{S}\times \wt{S})\mod(\DGCat).$$

\sssec{} \label{sss:Arinkin bis}

Note that since $S$ is smooth and $f$ is quasi-smooth, the scheme $\wt{S}$ is quasi-smooth, and so is 
$\wt{S}\underset{S}\times \wt{S}$. The localization
 $$\IndCoh_{\wt\CN}(\wt{S}\underset{S}\times \wt{S})\rightleftarrows \IndCoh(\wt{S}\underset{S}\times \wt{S})$$
 corresponds to a subset
 $$\wt\CN\subset \on{Sing}(\wt{S}\underset{S}\times \wt{S})$$
 (see \cite[Sect. 2.3]{AG} for what this means) defined as follows:
 
 \smallskip
 
 For a point $(\wt{s}_1,\wt{s}_2\,|\, f(\wt{s}_1)=s=f(\wt{s}_2)$) of $\wt{S}\underset{S}\times \wt{S}$, we have a naturally defined map
 \begin{multline*}
 \on{Sing}(\wt{S}\underset{S}\times \wt{S})_{(\wt{s}_1,\wt{s}_2)}= 
 H^{-1}(T^*_{(\wt{s}_1,\wt{s}_2)}(\wt{S}\underset{S}\times \wt{S}))\twoheadrightarrow \\
\to  \on{ker}(df:T^*_{s}(S)\to H^0(T_{\wt{s}_1}^*(\wt{S})))\cap  \on{ker}(df:T^*_{s}(S)\to H^0(T_{\wt{s}_2}^*(\wt{S})))
\hookrightarrow T^*_{s}(S).
\end{multline*}
 
 \smallskip
 
 We let $\wt\CN_{(\wt{s}_1,\wt{s}_2)}\subset \on{Sing}(\wt{S}\underset{S}\times \wt{S})_{(\wt{s}_1,\wt{s}_2)}$
 be the preimage of $\CN_s\subset T^*_{s}(S)$ under the above map. 
 
 \sssec{} \label{sss:obj in 2IndCoh}
 
Here is an example of an object in $\tIndCoh_\CN(S)$. Let $g:Z\to S$ be another quasi-smooth map. 
Set $\wt{Z}:=\wt{S}\underset{S}\times Z$. By the same recipe as in \secref{sss:Arinkin bis}, to $\CN$
we associate a subset
$$\wt\CN_Z\subset \on{Sing}(\wt{Z}).$$

Then $\IndCoh_{\wt\CN_Z}(\wt{Z})$ is naturally a module over $\IndCoh_{\wt\CN}(\wt{S}\underset{S}\times \wt{S})$.
Denote the resulting object of $\tIndCoh_\CN(S)$ by
$$\IndCoh_\CN(Z\mid S).$$

For example, if $\CN=\{0\}$, then $\IndCoh_\CN(Z\mid S)=\IndCoh(Z)$, viewed as a module over $\QCoh(S)$. 

\sssec{}

The point is that the above constructions are canonically independent of the choice of $\wt{S}$.
Namely, for two
$$\wt{S}_1\overset{f_1}\to S \overset{f_2}\leftarrow \wt{S}_2,$$
the bimodule given by
$$\IndCoh_{\wt\CN_{1,2}}(\wt{S}_1\underset{S}\times \wt{S}_2),$$
defines a Morita equivalence between 
$$\IndCoh_{\wt\CN}(\wt{S}_1\underset{S}\times \wt{S}_1) \text{ and } \IndCoh_{\wt\CN}(\wt{S}_2\underset{S}\times \wt{S}_2),$$ 
where $\wt\CN_{1,2}$ is the corresponding subset of $\on{Sing}(\wt{S}_1 \underset{S}\times\wt{S}_2)$. 

\sssec{} \label{sss:2-IndCoh stack}

By construction, $\tIndCoh(S)$ carries a monoidal action of the monoidal 2-category $\tQCoh(S)$. Using this action, the above construction
can be sheafified in the smooth topology, so $\tIndCoh_\CN(\CY)$ makes sense for a smooth stack $\CY$ and
$\CN$ a conical sub-Lagrangian subset of $T^*(\CY)$.

\smallskip

The construction in \secref{sss:obj in 2IndCoh} localizes as well: to a quasi-smooth map $\CZ\to \CY$ we can associate an object
$$\IndCoh_\CN(\CZ\mid\CY)\in \tIndCoh_\CN(\CY).$$

\smallskip

Furthermore, this construction make sense when $\CY$ is a formally smooth \emph{formal} algebraic stack. 

\sssec{}

We take $\CY:=\LS^{\on{restr}}_\cG(\ocD)$. Note that for $\sigma\in \LS^{\on{restr}}_\cG(\ocD)$, the cotangent fiber 
$T^*_\sigma(\LS^{\on{restr}}_\cG(\ocD))$ identifies with 
$$H^0(\ocD,\cg_\sigma).$$

Hence, we have a well-defined conical subset  
$$\Nilp\subset T^*(\LS^{\on{restr}}_\cG(\ocD)).$$

It is not difficult to show that $\Nilp$ is in fact Lagrangian. Thus, we can consider the 2-category $\tIndCoh_\Nilp(\LS^{\on{restr}}_\cG(\ocD))$.

\sssec{}

The following is the statement of the restricted local geometric Langands conjecture:

\begin{conj} \label{c:local GLC}
There exists an equivalence of 2-categories
$$\fL(G)\Cat_{\on{restr}}\simeq \tIndCoh_\Nilp(\LS^{\on{restr}}_\cG(\ocD)).$$
\end{conj} 

\smallskip

A very significant progress towards this conjecture has been recently achieved in the work of
G.~Dhillon, Y.~Varshavsky and D.~Yang, \cite{DVY} (taking the input from the depth-$0$ theory
developed in \cite{DLYZ}). 

\begin{rem}

As of now, the proof of \conjref{c:local GLC} proposed in \cite{DVY} has the following shortcomig:
it establishes an equivalence between the two sides as \emph{abtract} 2-categories. However,
both sides have additional structure:

\smallskip

As was mentioned above, $\tIndCoh_\Nilp(\LS^{\on{restr}}_\cG(\ocD))$ carries an action of
the monoidal 2-category $\tQCoh(\LS^{\on{restr}}_\cG(\ocD))$, i.e., $\QCoh(\LS^{\on{restr}}_\cG(\ocD))$ maps to the Bernstein center of 
$\tIndCoh_\Nilp(\LS^{\on{restr}}_\cG(\ocD))$.

\smallskip

Similarly, the \emph{fusion} construction defines a map from $\QCoh(\LS^{\on{restr}}_\cG(\ocD))$ to the 
Bernstein center of $\fL(G)\Cat_{\on{restr}}$.

\smallskip

What we do not know yet is why the equivalence of \cite{DVY} is compatible with the maps from 
$\QCoh(\LS^{\on{restr}}_\cG(\ocD))$ to the Bernstein centers of the two sides.

\end{rem}  

\sssec{Examples} \label{sss:examples local GLC}

Here are some examples of how the equivalence in \conjref{c:local GLC} works:

\smallskip

Consider the object $\ul\Shv(\Gr_G)\in \fL(G)\Cat_{\on{restr}}$. Then its image in $\tIndCoh_\Nilp(\LS^{\on{restr}}_\cG(\ocD))$
is $$\IndCoh_\Nilp(\LS^{\on{restr}}_\cG(\cD)\mid \LS^{\on{restr}}_\cG(\ocD)),$$ where we note that 
$\LS^{\on{restr}}_\cG(\cD)\simeq \on{pt}/\cG$. 

\smallskip 

Consider the object $\ul\Shv(\Fl^{\on{aff}}_G)\in \fL(G)\Cat_{\on{restr}}$. Then its image in $\tIndCoh_\Nilp(\LS^{\on{restr}}_\cG(\ocD))$
is $$\IndCoh_\Nilp\left(\LS^{\on{restr}}_\cT(\cD)\underset{\LS^{\on{restr}}_\cT(\ocD)}\times 
\LS^{\on{restr}}_\cB(\ocD)\mid \LS^{\on{restr}}_\cG(\ocD)\right).$$

\ssec{Consequences for the 1-categorical equivalence}

In this subsection we will perform a procedure parallel to one in \secref{ss:Tr of GLC},
but applied in the 2-categorical setting.

\sssec{}

Let us assume the equivalence of \conjref{c:local GLC}, i.e., the equivalence 
\begin{equation} \label{e:local GLC again}
\fL(G)\Cat_{\on{restr}}\simeq \tIndCoh_\Nilp(\LS^{\on{restr}}_\cG(\ocD)),
\end{equation}
as well as its functoriality with respect to $k$.

\smallskip

Taking $k=\ol\BF_q$, we obtain that both sides carry an action of Frobenius, and as a consequence we
obtain an equivalence 
\begin{equation} \label{e:Tr local GLC}
\Tr_{\DGCat}(\Frob,\fL(G)\Cat_{\on{restr}})\simeq \Tr_{\DGCat}(\Frob,\tIndCoh_\Nilp(\LS^{\on{restr}}_\cG(\ocD)))
\end{equation}
as objects of $\DGCat$, i.e., as DG categories. We will now calculate both sides of \eqref{e:Tr local GLC}.

\sssec{} \label{sss:Arinkin funct}

Let us return to the setting of Sects. \ref{sss:Arinkin}-\ref{sss:2-IndCoh stack}. Let
$$\phi:\CY_1\to \CY_2$$
be a map between smooth stacks. Let $\CN_i\subset T^*(\CY_i)$ be such that for 
$y_1\in \CY_1$, $y_2=\phi(y_1)\in \CY_2$ and the corresponding map
$d\phi:T^*_{y_2}(\CY_2)\to T^*_{y_1}(\CY_1)$, we have
$$d\phi^{-1}(\CN_{1,y_1})\subset  \CN_{2,y_2}.$$

\smallskip

We claim that in this case, we have a well-defined functor
\begin{equation} \label{e:pushforward tIndCoh}
\phi_*:\tIndCoh_{\CN_1}(\CY_1)\to \tIndCoh_{\CN_2}(\CY_2).
\end{equation} 

\sssec{}

The functor \eqref{e:pushforward tIndCoh} is obtained by smooth localization
of the following construction:

\smallskip

For $\wt{S}_2\to S_2$ as in \secref{sss:Arinkin bis}, consider the pullback square
$$
\CD
\wt{S}_1 @>{\wt\phi}>> \wt{S}_2 \\
@VVV @VVV \\
S_1 @>{\phi}>> S_2.
\endCD
$$

In this case, !-pullback gives rise to a functor
\begin{equation} \label{e:conf functor}
\IndCoh_{\wt\CN_2}(\wt{S}_2\underset{S}\times \wt{S}_2)\to \IndCoh_{\wt\CN_1}(\wt{S}_1\underset{S}\times \wt{S}_1)
\end{equation} 
that makes the diagram 
$$
\CD
\IndCoh(\wt{S}_2\underset{S}\times \wt{S}_2) @>>> \IndCoh(\wt{S}_1\underset{S}\times \wt{S}_1) \\
@VVV @VVV \\
\IndCoh_{\wt\CN_2}(\wt{S}_2\underset{S}\times \wt{S}_2) @>>> \IndCoh_{\wt\CN_1}(\wt{S}_1\underset{S}\times \wt{S}_1) 
\endCD
$$
commute\footnote{This follows from \cite[Proposition 7.1.3(b)]{AG}}. Restriction along \eqref{e:conf functor} defines the sought-for functor \eqref{e:pushforward tIndCoh}. 

\sssec{}

In the setting of \secref{sss:Arinkin funct} take $\CY_1=\CY_2=\CY$, so it makes sense to consider
\begin{equation} \label{e:Tr on 2-IndCoh}
\Tr_{\DGCat}(\phi_*,\tIndCoh_{\CN}(\CY))\in \DGCat
\end{equation}

\smallskip

The following is a variant of the calculation of \cite[Theorem 3.3.1]{BZNP}:

\begin{prop} \label{p:Tr on 2-IndCoh}
The object \eqref{e:Tr on 2-IndCoh} identifies with
$$\IndCoh_{\CN^\phi}(\CY^\phi),$$
where:

\begin{itemize}

\item $\CY^\phi$ is the fixed-point locus of $\phi$ (which is a quasi-smooth stack);

\item $\CN^\phi$ is the locus $\{(y\in \CY,\xi\in T^*_y(\CY))\,|\, \phi(y)=y,\, d\phi(\xi)=\xi\}$, which
is naturally a subset of $\on{Sing}(\CY^\phi)$.

\end{itemize}

\end{prop}

The above discussion is equally applicable when $\CY$ is a formally smooth formal algebraic stack.

\sssec{}

We apply \propref{p:Tr on 2-IndCoh} to $\CY=\LS^{\on{restr}}_\cG(\ocD)$ and $\phi=\Frob$. Denote
$$\LS^{\on{arthm}}_\cG(\ocD):=(\LS^{\on{restr}}_\cG(\ocD))^{\Frob}.$$

The subset $\Nilp\subset T^*(\LS^{\on{restr}}_\cG(\ocD))$ is Frobenius-invariant. Moreover, it is easy
to see that the embedding
$$\Nilp^{\Frob}\subset \on{Sing}(\LS^{\on{arthm}}_\cG(\ocD))$$
is an equality. Hence, from \propref{p:Tr on 2-IndCoh}, we obtain:

\begin{cor} \label{c:Tr on 2-IndCoh}
There exists a canonical isomorphism 
$$\Tr_{\DGCat}(\Frob,\tIndCoh_\Nilp(\LS^{\on{restr}}_\cG(\ocD))) \simeq \IndCoh(\LS^{\on{arthm}}_\cG(\ocD)).$$
\end{cor}

\sssec{}

Combing \eqref{e:Tr local GLC}, \corref{c:Tr on 2-IndCoh} and \conjref{c:local Tr conj}, we obtain:

\begin{corconj} \label{c:1-cat local Langlands}
There exists a canonical equivalence
$$\Shv(\on{Isoc}_G) \simeq \IndCoh(\LS^{\on{arthm}}_\cG(\ocD)).$$
\end{corconj}

We regard the statement of \corref{c:1-cat local Langlands} as an ultimate form of the local \emph{classical} 
Langlands conjecture. 

\begin{rem}

As was mentioned in the preamble to this section, \conjref{c:1-cat local Langlands} 
was first proposed by X.~Zhu in \cite{Zhu1} (see Conjecture 4.6.4 in {\it loc. cit.}).  The tame/depth-$0$
part of this conjecture has been recently established in \cite{Zhu2}. 

\end{rem}

\sssec{}

Recall now that in their work \cite{FS}, L.~Fargues and P.~Scholze suggested another form of the \emph{classical} 
Langlands conjecture. It reads as an equivalence
\begin{equation} \label{e:FS}
D(\Bun^{\on{local}}_G)\simeq \IndCoh(\on{Param}_\cG),
\end{equation} 
where:

\begin{itemize}

\item $\Bun^{\on{local}}_G$ is the analytic stack of $G$-bundles on the Fargues-Fontaine curve, and $D(-)$
refers to an appropriately defined category of $\ol\BQ_\ell$-adic sheaves;

\item $\on{Param}_\cG$ is the stack of Langlands parameters, as defined in \cite[Sect. I.8]{FS}.

\end{itemize}

However:

\begin{itemize}

\item One can show that the stacks  $\on{Param}_\cG$ is canonically isomorphic to $\LS^{\on{arthm}}_\cG(\ocD)$;

\item The ongoing work \cite{GIJHZ} aims to establish the equivalence $D(\Bun^{\on{local}}_G)\simeq \Shv(\on{Isoc}_G)$.

\end{itemize}

Thus, the statement of \eqref{e:FS} is equivalent to that of \corref{c:1-cat local Langlands}. 

\section{The global ramified case}  \label{s:global}

The contents of this section constitute work-in-progress, joint with S.~Raskin and Y.~Varshavsky. 

\smallskip

We will propose a set of conjectures about the global ramified Langlands theory, assuming
the local Langlands theory discussed in the previous section. There are three main conjectures:

\begin{itemize}

\item A spectral description of the global automorphic category, \conjref{c:global ramified} (this is a 
ramified counterpart of \conjref{c:restr GLC});

\item A global trace isomorphism, \conjref{c:Tr isom ramified} (this is a ramified counterpart of
\thmref{t:Tr isom});

\item A spectral description of the (enhanced version of) the space of automorphic functions, \conjref{c:cl in 2-IndCoh}
(this is a ramified counterpart of \corref{c:Autom descr}). 

\end{itemize}

The last conjecture is a formal consequence of the first two and supplies an 
answer to the basic question of the Langlands program. 

\ssec{The global category in the ramified case} \label{ss:HL ramified}

In this subsection we introduce the key geometric player: the category $\Shv_{\on{HL}}(\Bun_G^{\on{level}_{\ul{x}}})$,
or, rather, its version $\ul\Shv_{\on{HL}}(\Bun_G^{\on{level}_{\ul{x}}})\in \AGCat$. 

\sssec{}

Let $X$ be again a smooth and complete curve over the ground field $k$. We fix a finite set of
$k$-points $\ul{x}\subset X$. Denote $\oX:=X-\ul{x}$.

\smallskip

Let $\Bun_G^{\on{level}_{\ul{x}}}$ denote the moduli stack of $G$-bundles
on $X$ with a full level structure at $\ul{x}$. If $\ul{x}\neq \emptyset$, then $\Bun_G^{\on{level}_{\ul{x}}}$
is actually a scheme (of infinite type). 

\smallskip

The stack $\Bun_G^{\on{level}_{\ul{x}}}$ carries a natural action of $\fL(G)_{\ul{x}}$, i.e., the product of
copies of $\fL(G)$ over the points that comprise $\ul{x}$.

\smallskip

Our initial object of study is the DG category $\Shv(\Bun_G^{\on{level}_{\ul{x}}})$, or rather 
\begin{equation} \label{e:Shv BunG level AG}
\ul\Shv(\Bun_G^{\on{level}_{\ul{x}}})\in \AGCat.
\end{equation} 

The latter carries a natural action of $\fL(G)_{\ul{x}}$, i.e., can be thought of as object of
$\fL(G)_{\ul{x}}\Cat$.

\sssec{}

Let us temporarily place ourselves in the context of D-modules, in which case $\Dmod(\Bun_G^{\on{level}_{\ul{x}}})$
is naturally an object of $\fL(G)_{\ul{x}}\Cat$. 

\smallskip

Unfortunately, with the current state of knowledge, there is very little that we can say about this category (except
when $G=T$ is a torus). For example, we do not know how to answer such basic questions as to whether or not
$\Dmod(\Bun_G^{\on{level}_{\ul{x}}})$ is compactly generated or dualizable. 

\sssec{}

The object \eqref{e:Shv BunG level AG} carries an action of the Hecke functors
$$\on{H}_V:\ul\Shv(\Bun_G^{\on{level}_{\ul{x}}})\to \ul\Shv(\Bun_G^{\on{level}_{\ul{x}}}\times X),$$
i.e., the functors
$$\on{H}_V:\Shv(\Bun_G^{\on{level}_{\ul{x}}}\times S)\to \Shv(\Bun_G^{\on{level}_{\ul{x}}}\times S\times \oX).$$

Let 
$$\ul\Shv_{\on{HL}}(\Bun_G^{\on{level}_{\ul{x}}})\in \AGCat$$ 
be the object, whose value on $S\in \Sch^{\on{aff}}$ is the full subcategory of $\Shv(\Bun_G^{\on{level}_{\ul{x}}}\times S)$ consisting 
of objects $\CF$ for which
$$\on{H}_V(\CF)\in  \Shv(\Bun_G^{\on{level}_{\ul{x}}}\times S)\otimes \qLisse(\oX)\subset  \Shv(\Bun_G^{\on{level}_{\ul{x}}}\times S\times \oX).$$

\smallskip

Denote
$$\Shv_{\on{HL}}(\Bun_G^{\on{level}_{\ul{x}}}):=\ul\Shv_{\on{HL}}(\Bun_G^{\on{level}_{\ul{x}}})_{\on{pt}}.$$

\sssec{}

It is clear that $\ul\Shv_{\on{HL}}(\Bun_G^{\on{level}_{\ul{x}}})$ inherits an action of $\fL(G)_{\ul{x}}$, i.e., 
lifts to an object of $\fL(G)_{\ul{x}}\Cat$.

\smallskip

The following is one of our key conjectures:

\begin{conj} \label{c:global restr}
The object $\ul\Shv_{\on{HL}}(\Bun_G^{\on{level}_{\ul{x}}})\in \fL(G)_{\ul{x}}\Cat$ belongs to $\fL(G)_{\ul{x}}\Cat_{\on{restr}}$.
\end{conj} 

\sssec{} \label{sss:check global restr}

Let us run a reality check on \conjref{c:global restr}. Consider the objects
$$\ul\Shv_{\on{HL}}(\Bun_G^{\on{level}_{\ul{x}}}) \text{ and } \ul\Shv(\Gr_{G,\ul{x}}) \text{ in } \fL(G)_{\ul{x}}\Cat.$$
Since both are supposed to be restricted, by \eqref{e:inv restr}, the object
$$\ul\Shv_{\on{HL}}(\Bun_G) \simeq \left(\ul\Shv_{\on{HL}}(\Bun_G^{\on{level}_{\ul{x}}}) \otimes \ul\Shv(\Gr_{G,\ul{x}}) \right)^{\fL(G)_{\ul{x}}}\in \AGCat$$
is supposed belongs to the essential image of $\DGCat$ under $\on{emb.restr}$. However, the latter is indeed the case by \cite[Theorem 1.3.7]{AGKRRV2}.

\sssec{}

We also add the following technical conjecture:

\begin{conj} \label{c:HL is dualizable}
The object $\ul\Shv_{\on{HL}}(\Bun_G^{\on{level}_{\ul{x}}})$ is dualizable in $\fL(G)_{\ul{x}}\Cat$\footnote{Or, 
or, equivalently, in $\AGCat$.}.
\end{conj} 

Note that if we assume \conjref{c:global restr}, we can restate \conjref{c:HL is dualizable} as saying that
$\ul\Shv_{\on{HL}}(\Bun_G^{\on{level}_{\ul{x}}})$ is dualizable as an object of $\fL(G)_{\ul{x}}\Cat_{\on{restr}}$.

\ssec{An enhancement of the space of automorphic functions}

The traditional object of study in the theory of automorphic functions is .... \emph{the space of automorphic functions},
which is a representation of the ad\`ele group, or, for a fixed locus of ramification, an object of
$\Rep(\fL(G)_{\ul{x}}(\BF_q))$.  

\smallskip

However, as we have seen in the previous section, the latter category is not quite
the right object\footnote{E.g., if we want to match it with the spectral side.}, and should be enhanced to $\Shv(\on{Isoc}_G)$.
In this subsection we construct the corresponding enlargement of the space of automorphic functions. 

\sssec{} \label{sss:autom enh}

Let $k=\ol\BF_q$, but assume that $G$, $X$ and $\ul{x}$ are defined over $\BF_q$. Hence, it makes sense to consider
the set $\Bun_G^{\on{level}_{\ul{x}}}(\BF_q)$ (of isomorphism classes of) $\BF_q$-points of $\Bun_G^{\on{level}_{\ul{x}}}$.

\smallskip

The stack $\Bun_G^{\on{level}_{\ul{x}}}$ is naturally the inverse limit of a family of stacks locally of finite type. Hence,
 $\Bun_G^{\on{level}_{\ul{x}}}(\BF_q)$ carries a natural topology. We let
 \begin{equation} \label{e:ramified autom functions}
 \on{Autom}(X,G)^{\on{level}_{\ul{x}}}
 \end{equation} 
 denote the space of compactly supported locally constant functions on $\Bun_G^{\on{level}_{\ul{x}}}(\BF_q)$.
 This is the space of (compactly supported, smooth) automorphic functions of full level at $\ul{x}$.
 
 \smallskip
 
 The vector space $\on{Autom}(X,G)^{\on{level}_{\ul{x}}}$ carries a smooth action of the group $\fL(G)_{\ul{x}}(\BF_q)$. 
  
 \sssec{}
 
 We will now define an object
 \begin{equation} \label{e:Autom enh}
 \on{Autom}(X,G)^{\on{enh}_{\ul{x}}}\in \Shv(\on{Isoc}_{G,\ul{x}}),
 \end{equation} 
 from which \eqref{e:ramified autom functions} can be recovered by applying the functor $\iota^!$ of
 \eqref{e:unit stratum}. In the above formula, the space $\on{Isoc}_{G,\ul{x}}$ is the product of copies
 of $\on{Isoc}_G$ over points that comprise $\ul{x}$. 
 
 \smallskip
 
 Recall (see \secref{sss:Isoc geom}) that the category $\Shv(\on{Isoc}_{G,\ul{x}})$ is self-dual, so the datum of an object 
 in it is equivalent to that of a colimit-preserving functor
 \begin{equation} \label{e:functor out of Isoc}
 \Shv(\on{Isoc}_{G,\ul{x}})\to \Vect.
 \end{equation} 
 
 So, our goal is to define such a functor that would correspond to the sough-for object \eqref{e:Autom enh}.
 
 \sssec{}
 
 Let $\Sht_{\ul{x}}$ denote the stack of \emph{shtukas with a leg at $\ul{x}$}, i.e., 
 $$\Sht_{\ul{x}}:=\Bun_G\underset{\Bun_G\times \Bun_G}\times \CH_{\ul{x}}.$$
 where:
 
\begin{itemize}
 
\item The map $\Bun_G\to \Bun_G\times \Bun_G$ has components $(\on{Id},\Frob)$;

\item $\CH_{\ul{x}}$ is the Hecke groupoid at $x$, i.e., the fiber product
$$\Bun_G\underset{\Bun_G(\oX)}\times \Bun_G.$$

 \end{itemize}
 
 \smallskip
 
 The stack $\Sht_{\ul{x}}$ carries a canonical action of $\CH_{\ul{x}}$; denote 
 $$\oSht_{\ul{x}}:=\Sht_{\ul{x}}/\CH_{\ul{x}}.$$
 
 Note that since $\CH_{\ul{x}}$ is proper, the category $\Shv(\oSht_{\ul{x}})$ is well-defined.
 Moreover, we have a well-defined functor 
 $$\on{C}^\cdot_c(\oSht_{\ul{x}},-):\Shv(\oSht_{\ul{x}})\to \Vect.$$
 
 \begin{rem}
 
 We can think of $\oSht_{\ul{x}}$ as $(\Bun_G(\oX))^{\Frob}$. 
 
 \end{rem}
 
 \sssec{}
 
 Restriction from $X$ to the multi-disc around $\ul{x}$ gives rise to a map
 $$\pi:\oSht_{\ul{x}}\to \on{Isoc}_{G,\ul{x}}.$$
 
We define the functor \eqref{e:functor out of Isoc} by sending $\CF\in \Shv(\on{Isoc}_{G,\ul{x}})$ to
$$\on{C}^\cdot_c(\oSht_{\ul{x}},\pi^!(\CF)).$$
 
 \sssec{}
 
 Let us show how to construct an identification
 \begin{equation} \label{e:fiber of enh}
 \iota^!(\on{Autom}(X,G)^{\on{enh}_{\ul{x}}})\simeq \on{Autom}(X,G)^{\on{level}_{\ul{x}}}.
 \end{equation} 
 
Let $\delta_1$ be the delta-function object in $\Shv(\on{Isoc}_{G,\ul{x}})$ corresponding to the
unit point of $\on{Isoc}_{G,\ul{x}}$. I.e.,
$$\delta_1:=\iota_!(R_{\fL(G)_{\ul{x}}(\BF_q)}),$$
where $R_{\fL(G)_{\ul{x}}(\BF_q)}\in \Rep(\fL(G)_{\ul{x}}(\BF_q))$ is regular representation, i.e., the space of
smooth compactly supported functions on $\fL(G)_{\ul{x}}(\BF_q)$. Note that by construction, $\delta_1$ 
is equipped with a natural action of $\fL(G)_{\ul{x}}(\BF_q)$.
 
 \smallskip
 
 Now, it follows from the definition of the self-duality on $\Shv(\on{Isoc}_{G,\ul{x}})$ that the functor $\iota^!$
 identifies with the functor of pairing with  $\delta_1$. Thus, in order to prove \eqref{e:fiber of enh}, we have to
 establish an isomorphism
 $$\on{C}^\cdot_c(\oSht_{\ul{x}},\pi^!(\delta_1))\simeq \on{Autom}(X,G)^{\on{level}_{\ul{x}}}.$$
 
 This follows by base change from the Cartesian diagram\footnote{The horizontal arrows in this diagram are ind-pro-proper, so base change holds.}
 $$
 \CD
 (\Bun_G^{\on{level}_{\ul{x}}})^{\Frob} @>>> \oSht_{\ul{x}} \\
 @VVV @VVV \\
 \on{pt} @>>> \on{Isoc}_{G,\ul{x}}.
 \endCD
 $$
 
 \ssec{The trace conjecture in the ramified case}
 
In this subsection we we retain the setting of \secref{sss:autom enh}. 
We will state a conjecture that generalizes \thmref{t:Tr isom} in the 
 ramified situation.
 
 \sssec{}
 
 Consider the object 
 $$\ul\Shv_{\on{HL}}(\Bun_G^{\on{level}_{\ul{x}}})\in \fL(G)_{\ul{x}}\Cat.$$
 
 It carries a natural action of Frobenius, which makes it equivariant with respect to the action of $\Frob$
 on $\fL(G)_{\ul{x}}\Cat$. Moreover, according to \conjref{c:HL is dualizable}, it is dualizable, so it makes
 sense to consider
 $$\on{cl}(\ul\Shv_{\on{HL}}(\Bun_G^{\on{level}_{\ul{x}}}),\Frob)\in \Tr_{\AGCat}(\Frob,\fL(G)_{\ul{x}}\Cat)_{\on{pt}}\simeq
 \Shv(\on{Isoc}_G).$$
 
 Note that if we assume \conjref{c:global restr}, we can interpret the latter object also as a class taken in
 $$\Tr_{\DGCat}(\Frob,\fL(G)_{\ul{x}}\Cat_{\on{restr}})\overset{\text{\conjref{c:local Tr conj}}}\simeq  \Shv(\on{Isoc}_G).$$

\sssec{}

We propose:

\begin{conj} \label{c:Tr isom ramified}
There exists a canonical isomorphism in $\Shv(\on{Isoc}_G)$:
$$ \on{Autom}(X,G)^{\on{enh}_{\ul{x}}}\simeq \on{cl}(\ul\Shv_{\on{HL}}(\Bun_G^{\on{level}_{\ul{x}}}),\Frob).$$
\end{conj}

\begin{rem}

It is a somewhat non-trivial exercise (one that uses \eqref{e:inv restr}) to show that 
that the isomorphism of \conjref{c:Tr isom ramified} implies the isomorphism
of \thmref{t:Tr isom}.

\end{rem} 

\ssec{The geometric Langlands equivalence in the global ramified case}

\sssec{} \label{sss:geom Langlands ram intro}

We return to the set-up of \secref{ss:HL ramified}. The goal of the (restricted) version of the global 
geometric Langlands conjecture in the presence of ramification is to describe the object $\ul\Shv_{\on{HL}}(\Bun_G^{\on{level}_{\ul{x}}})$ 
in spectral terms.

\smallskip

However, we cannot quite do it ``as-is": first, $\ul\Shv_{\on{HL}}(\Bun_G^{\on{level}_{\ul{x}}})$ is not
even a category (but rather an object of $\AGCat$); and second, it carries an action of
$\fL(G)_{\ul{x}}$, and objects arising from $\cG$ do not have room for such a structure.

\smallskip

But this is a familiar phenomenon in the usual Langlands theory: the global ramified conjecture
is formulated in terms of the local conjecture, i.e., our aim should rather be to find the image
of $\ul\Shv_{\on{HL}}(\Bun_G^{\on{level}_{\ul{x}}})$ under the local Langlands equivalence.

\sssec{}

Let us assume the validity of \conjref{c:global restr}, so we can view 
\begin{equation} \label{e:HL in restr}
\ul\Shv_{\on{HL}}(\Bun_G^{\on{level}_{\ul{x}}})\in \fL(G)_{\ul{x}}\Cat_{\on{restr}}.
\end{equation}

Recall now that, according to \conjref{c:local GLC}, we have an equivalence of 2-categories:
\begin{equation} \label{e:local GLC again x}
\fL(G)_{\ul{x}}\Cat_{\on{restr}}\simeq \tIndCoh_\Nilp(\LS^{\on{restr}}_\cG(\ocD_{\ul{x}})),
\end{equation}
where $\ocD_{\ul{x}}$ is the union of formal punctured discs around the points that comprise $\ul{x}$.

\smallskip

The global ramified geometric Langlands conjecture will amount to describing 
the image of \eqref{e:HL in restr} under \eqref{e:local GLC again x}.

\sssec{}

Consider the (formal) stack $\LS_\cG^{\on{restr}}(\oX)$ and the morphsim
\begin{equation} \label{e:from global to local LS}
\fr:\LS_\cG^{\on{restr}}(\oX)\to \LS_\cG^{\on{restr}}(\ocD_{\ul{x}}),
\end{equation} 
given by restriction. 

\smallskip

Consider the object
\begin{equation} \label{e:rel IndCoh}
\IndCoh_\Nilp(\LS_\cG^{\on{restr}}(\oX)\mid  \LS_\cG^{\on{restr}}(\ocD_{\ul{x}}))\in \tIndCoh_\Nilp(\LS^{\on{restr}}_\cG(\ocD_{\ul{x}})),
\end{equation} 
see \secref{sss:2-IndCoh stack} for the notation. 

\begin{conj} \label{c:global ramified}
Under the (conjectural) equivalence \eqref{e:local GLC again x}, the object \eqref{e:HL in restr} corresponds to \eqref{e:rel IndCoh}.
\end{conj} 

\sssec{Example}

Let us explain how \conjref{c:global ramified} recovers the unramified version, given by \conjref{c:restr GLC}:

\smallskip

We consider the following pairs of objects
$$\ul\Shv(\Gr_{G,\on{x}}) \text{ and } \ul\Shv_{\on{HL}}(\Bun_G^{\on{level}_{\ul{x}}})$$ 
in $\fL(G)_{\ul{x}}\Cat_{\on{restr}}$
and 
$$\IndCoh_\Nilp(\LS^{\on{restr}}_\cG(\cD)\mid \LS^{\on{restr}}_\cG(\ocD_{\ul{x}})),
\IndCoh_\Nilp(\LS_\cG^{\on{restr}}(\oX)\mid  \LS_\cG^{\on{restr}}(\ocD_{\ul{x}}))$$ in 
$\tIndCoh_\Nilp(\LS^{\on{restr}}_\cG(\ocD_{\ul{x}}))$.

\smallskip

They are supposed to correspond to one another under \eqref{e:local GLC again x}. Hence, we obtain an equivalence of DG categories
\begin{multline} \label{e:Hom from Sph}
\ul\Hom
%_{\fL(G)_{\ul{x}}\Cat_{\on{restr}}}
\left(\ul\Shv(\Gr_{G,\on{x}}), \ul\Shv_{\on{HL}}(\Bun_G^{\on{level}_{\ul{x}}})\right) \simeq \\
\simeq \ul\Hom
%_{\tIndCoh_\Nilp(\LS^{\on{restr}}_\cG(\ocD_{\ul{x}}))}
\left(\IndCoh_\Nilp(\LS_\cG^{\on{restr}}(\oX)\mid  \LS_\cG^{\on{restr}}(\ocD_{\ul{x}})),
\IndCoh_\Nilp(\LS_\cG^{\on{restr}}(\oX)\mid  \LS_\cG^{\on{restr}}(\ocD_{\ul{x}}))\right),
\end{multline}
where $\ul\Hom(-,-)$ stands for the enriched $\Hom$ over $\DGCat$, taken in 
$$\fL(G)_{\ul{x}}\Cat_{\on{restr}} \text{ and } \tIndCoh_\Nilp(\LS^{\on{restr}}_\cG(\ocD_{\ul{x}})),$$
respectively.

\smallskip

The left-hand side in \eqref{e:Hom from Sph} is isomorphic to
$$\ul\Shv_{\on{HL}}(\Bun_G)_{\on{pt}}= \Shv_{\on{HL}}(\Bun_G).$$

A direct calculation shows that the right-hand side in \eqref{e:Hom from Sph} is isomorphic to
$$\IndCoh_\Nilp(\LS_\cG(X)).$$

Thus, \eqref{e:Hom from Sph} gives us back the equivalence of \conjref{c:restr GLC}.

\ssec{Arithmetic consequences of the Langlands conjecture}

In this subsection we return to the setting of \secref{sss:autom enh}. We will 
combine all the ingredients developed so far, and try answer the following question:
what is the spectral description of the space $\on{Autom}(X,G)^{\on{level}_{\ul{x}}}$?

\sssec{}

As in \secref{sss:geom Langlands ram intro}, we cannot relate $\on{Autom}(X,G)^{\on{level}_{\ul{x}}}$, viewed
as a plain vector space, to the spectral side. Rather, we need to view it as an object of $\Rep(\fL(G)_{\ul{x}}(\BF_q))$. 

\smallskip

However, (un)surprisingly, the latter turns out to be not quite the right question (at least, not one that has a 
clean answer). The reason is that we do not know how to translate the category $\Rep(\fL(G)_{\ul{x}}(\BF_q))$.

\smallskip

What we \emph{will} be able to do is describe the object
$$\on{Autom}(X,G)^{\on{enh}_{\ul{x}}}\in \Shv(\on{Isoc}_G)$$
in terms of the equivalence 
\begin{equation} \label{e:1-cat local Langlands}
\Shv(\on{Isoc}_{G,\ul{x}}) \simeq \IndCoh(\LS^{\on{arthm}}_\cG(\ocD_{\ul{x}})))
\end{equation}
of \conjref{c:1-cat local Langlands}. 

\sssec{}

We return to the setting of \conjref{c:global ramified}, and we assume that the two objects
match in a way compatible with their Frobenius-equivariant structures.

\smallskip

Thus, we have an equivalence of pairs
\begin{multline} \label{e:pairs}
\ul\Shv_{\on{HL}}(\Bun_G^{\on{level}_{\ul{x}}})\in \fL(G)_{\ul{x}}\Cat_{\on{restr}}\simeq \\
\simeq \IndCoh_\Nilp(\LS_\cG^{\on{restr}}(\oX)\mid  \LS_\cG^{\on{restr}}(\ocD_{\ul{x}}))\in \tIndCoh_\Nilp(\LS^{\on{restr}}_\cG(\ocD_{\ul{x}})),
\end{multline} 
from which we obtain an equivalence 
\begin{multline} \label{e:pairs bis}
\on{cl}(\ul\Shv_{\on{HL}}(\Bun_G^{\on{level}_{\ul{x}}}),\Frob)\in \Tr_{\DGCat}(\Frob,\fL(G)_{\ul{x}}\Cat_{\on{restr}}) \simeq \\
\simeq \on{cl}\left(\IndCoh_\Nilp(\LS_\cG^{\on{restr}}(\oX)\mid  \LS_\cG^{\on{restr}}(\ocD_{\ul{x}})),\Frob\right) \in
\Tr(\Frob,\IndCoh_\Nilp(\LS_\cG^{\on{restr}}(\oX)\mid  \LS_\cG^{\on{restr}}(\ocD_{\ul{x}}))).
\end{multline}

In other words, the objects
\begin{equation} \label{e:geom pair}
\on{cl}(\ul\Shv_{\on{HL}}(\Bun_G^{\on{level}_{\ul{x}}}),\Frob)\in \Tr_{\DGCat}(\Frob,\fL(G)_{\ul{x}}\Cat_{\on{restr}}) 
\overset{\text{\eqref{e:local trace conj}}} \simeq \Shv(\on{Isoc}_{G,\ul{x}})
\end{equation}
\begin{multline}  \label{e:spectral pair}
\on{cl}\left(\IndCoh_\Nilp(\LS_\cG^{\on{restr}}(\oX)\mid  \LS_\cG^{\on{restr}}(\ocD_{\ul{x}})),\Frob\right) \in \\
\Tr(\Frob,\IndCoh_\Nilp(\LS_\cG^{\on{restr}}(\oX)\mid  \LS_\cG^{\on{restr}}(\ocD_{\ul{x}}))) \overset{\text{\corref{c:Tr on 2-IndCoh}}}\simeq
\IndCoh(\LS^{\on{arthm}}_\cG(\ocD_{\ul{x}}))
\end{multline}
correspond to one another under the equivalence \eqref{e:1-cat local Langlands}. 

\sssec{}

Let us identify explicitly the two objects \eqref{e:geom pair} and \eqref{e:spectral pair}. 

\smallskip

First, by \conjref{c:Tr isom ramified}, the object \eqref{e:geom pair} identifies with $\on{Autom}(X,G)^{\on{enh}_{\ul{x}}}$. 

\sssec{}

In order to identify \eqref{e:spectral pair}, we revisit the equivalence of \propref{p:Tr on 2-IndCoh}. Let
$g:\CZ\to \CY$ be a quasi-smooth map, and assume that the endomorphism $\phi$ of $\CY$ lifts to
an endomorphism of $\CZ$ (we denote it by the same symbol $\phi$). 

\smallskip

In this case, the operation of !-pullback defines a map
\begin{equation} \label{e:pushforward tIndCoh equiv}
\IndCoh_\CN(\CZ\mid \CY) \to \phi_*(\IndCoh_\CN(\CZ\mid \CY)),
\end{equation} 
where $\phi_*$ denotes the functor \eqref{e:pushforward tIndCoh}. 

\smallskip

We denote the map in \eqref{e:pushforward tIndCoh equiv} by $(\phi^!\mid\CY)$. Thus, we can consider the object
\begin{equation} \label{e:cl in 2-IndCoh}
\on{cl}(\IndCoh_\CN(\CZ\mid \CY),(\phi^!\mid\CY))\in \Tr(\phi_*,\tIndCoh_\CN(\CY)).
\end{equation}

The following is plausible (and, presumably, not difficult):

\begin{conj} \label{c:cl in 2-IndCoh}
Under the identification 
$$\Tr_{\DGCat}(\phi_*,\tIndCoh_{\CN}(\CY))\simeq \IndCoh_{\CN^\phi}(\CY^\phi)$$
of \propref{p:Tr on 2-IndCoh}, the object \eqref{e:cl in 2-IndCoh} corresponds to
$$((g^\phi)_*^\IndCoh(\omega_{\CZ^\phi}))_{\CN^\phi}\in \IndCoh_{\CN^\phi}(\CY^\phi),$$
where:

\begin{itemize}

\item $g^\phi$ denotes the morphism $\CZ^\phi\to \CY^\phi$;

\item $(-)_{\CN^\phi}$ denotes the image of a given object under the right adjoint to the embedding
$$\IndCoh_{\CN^\phi}(\CY^\phi)\hookrightarrow \IndCoh(\CY^\phi).$$

\end{itemize}

\end{conj} 

\sssec{}

We apply \conjref{c:cl in 2-IndCoh} to the morphism \eqref{e:from global to local LS} and $\phi=\Frob$. 
Denote
$$\LS^{\on{arithm}}_\cG(\oX):=(\LS^{\on{restr}}_\cG(\oX))^{\Frob}.$$

Denote by $\fr^{\on{arithm}}$ the resulting morphism
$\LS^{\on{arithm}}_\cG(\oX)\to \LS^{\on{arithm}}_\cG(\ocD_{\ul{x}})$.
Thus, from \conjref{c:cl in 2-IndCoh}, we obtain:

\begin{corconj} \label{cor:cl in 2-IndCoh}
The object \eqref{e:spectral pair} identifies with 
$$(\fr^{\on{arithm}})_*^\IndCoh(\omega_{\LS^{\on{arithm}}_\cG(\oX)})\in \IndCoh(\LS^{\on{arthm}}_\cG(\ocD_{\ul{x}})).$$
\end{corconj}

\sssec{}

Thus, combining, we obtain\footnote{The statement of \corref{c:ultimate} was inspired by a seminar talk by X.~Zhu in 2020.}:

\begin{corconj} \label{c:ultimate}
Under the equivalence \eqref{e:1-cat local Langlands}, the object 
$$\on{Autom}(X,G)^{\on{enh}_{\ul{x}}}\in \Shv(\on{Isoc}_G)$$
corresponds to the object
$$(\fr^{\on{arithm}})_*^\IndCoh(\omega_{\LS^{\on{arithm}}_\cG(\oX)})\in \IndCoh(\LS^{\on{arthm}}_\cG(\ocD_{\ul{x}})).$$
\end{corconj}

We regard the statement of \corref{c:ultimate} as an answer to the question of the spectral description of the 
space of automorphic functions in the ramified situation over function fields. 

\begin{rem}

Note that, according to \cite{Laf}, the object $\on{Autom}(X,G)^{\on{level}_{\ul{x}}}\in \Rep(\fL(G)_{\ul{x}}(\BF_q))$ carries am action of the \emph{excursion} algebra,
which we can interpret as $\Gamma(\LS^{\on{arithm}}_\cG(\oX),\CO_{\LS^{\on{arithm}}_\cG(\oX)})$. In a work-in-progress with 
A.~Genestier, A.~Eteve and V.~Lafforgue, we extend this action to one on $\on{Autom}(X,G)^{\on{enh}_{\ul{x}}}$ as an object
of $\Shv(\on{Isoc}_G)$. Under the identification of \corref{c:ultimate}, this action is supposed to correspond to the natural action of 
$\Gamma(\LS^{\on{arithm}}_\cG(\oX),\CO_{\LS^{\on{arithm}}_\cG(\oX)})$ on $\omega_{\LS^{\on{arithm}}_\cG(\oX)}$. 

\end{rem}

\appendix

\section{The trace formalism}

\ssec{Abstract framework for the trace} \label{ss:Tr abs}

Let $\bO$ be a symmetric monoidal category; let $\one_\bO\in \bO$ denote the unit object. 
:et $\bo\in \bO$ be a dualizable object, i.e., there exists an object $\bo^\vee$, equipped
with the maps
$$\one_\bO\overset{\on{unit}}\to \bo\otimes \bo^\vee \text{ and } \bo\otimes \bo^\vee\overset{\on{counit}}\to \one_\bO,$$
satisfying the usual axioms. 

\smallskip

Let $F$ be an endomorphism of $\bo$. The trace $\Tr(F,\bo)\in \End_\bO(\one_\bO)$ is defined as the composition
$$\one_{\bO}\overset{\on{unit}}\to \bo\otimes \bo^\vee\overset{F\otimes \on{id}}\to \bo\otimes \bo^\vee \overset{\on{counit}}\to \one_\bO.$$

\ssec{The trace of an endofunctor on a category} \label{ss:cat Tr}

Our main example will be $\bO=\DGCat$ (see \cite[Chapter 1, Sect. 10]{GaRo}. The unit is the DG category of vector spaces,
denoted $\Vect$. 

\smallskip

Thus, it makes sense to talk about the trace of an endofunctor $F$ of a dualizable\footnote{The higher we go in the
categorical hierarchy, the easier it becomes for objects to be dualizable. E.g., a vector space is dualizable if and only if it is finite-dimensional.
By contrast, most categories that we encounter in practice are dualizable.} DG category $\bC$
$$\Tr(F,\bC)\in \Vect.$$

\ssec{Trace on $\QCoh$} \label{sss:Tr QCoh}

Here is a typical trace calculation in $\DGCat$. Let $\CY$ be a prestack with a schematic diagonal, and such that the functor
$$\QCoh(\CY)\otimes \QCoh(\CY)\overset{\boxtimes}\to \QCoh(\CY\times \CY)$$
is an equivalence. Then $\QCoh(\CY)$ is dualizable and self-dual with the unit being
$$\Vect\simeq \QCoh(\on{pt}) \overset{k\mapsto CO_\CY}\longrightarrow \QCoh(\CY)\overset{(\Delta_\CY)_*}\longrightarrow \QCoh(\CY\times \CY)\simeq \QCoh(\CY)\otimes \QCoh(\CY)$$
and the counit
$$\QCoh(\CY)\otimes \QCoh(\CY)\overset{\sim}\to \QCoh(\CY\times \CY)\overset{(\Delta_\CY)^*}\longrightarrow \QCoh(\CY)\overset{\Gamma(\CY,-)}\to \Vect.$$

Let $\phi$ be an endomorphism of $\CY$ (assumed schematic as a map) and consider the endofunctor of $\QCoh(\CY)$ given by $\phi_*$. 
We claim that
$$\Tr(\phi_*,\QCoh(\CY))\simeq \Gamma(\CY^\phi,\CO_{\CY^\phi}),$$
where $\CY^\phi$ is the fixed-point locus of $\phi$ on $\CY$. Indeed, this follows by base change from the following diagram
$$
\CD
\CY^\phi & @>>> & \CY @>>> \on{pt} \\
@VVV & & @VV{\Delta_\CY}V   \\
\CY @>{\Delta_\CY}>> \CY\times \CY @>{\phi\times \on{id}}>> \CY\times \CY \\
@VVV \\
\on{pt},
\endCD
$$
in which we take pullbacks with respect to vertical arrows and pushforwards with respect to horizontal arrows. 

\ssec{Functoriality of trace} \label{sss:funct Tr}

Note that $\DGCat$ has an additional structure, namely, that of a 2-category. This structure leaves room for
additional functoriality of the trace construction. 

\smallskip

Let us return to the setting of \secref{ss:Tr abs}, where $\bO$ is now a symmetric monoidal 2-category. 
Let $\bo_1,\bo_2$ be a pair of dualizable objects, each equipped with an endomorphism $F_i$, and let
$\Phi:\bo_1\to \bo_2$ be a map that \emph{weakly} commutes with the endomorphisms, i.e., we are given
a natural transformation
$$\Phi\circ F_1\overset{\alpha}\to F_2\circ \Phi.$$

Assume now that $\Phi$ admits a \emph{right adjoint} in $\bO$. Then we claim that there exists a naturally defined map
$$\Tr(F_1,\bo_1)\to \Tr(F_2,\bo_2)$$
in $\End(\one_\bO)$. Namely, this map equals the composition
$$\Tr(F_1,\bo_1) \overset{\on{unit}}\to \Tr(\Phi^R\circ \Phi \circ F_1,\bo_1) \overset{\on{cyclicity}}\simeq 
\Tr(\Phi \circ F_1\circ \Phi^R,\bo_2)\overset{\alpha}\to \Tr(F_2\circ \Phi \circ \Phi^R,\bo_2) \overset{\on{counit}}\to \Tr(F_2,\bo_2).$$

\ssec{The class of an object} \label{sss:class}

A particular case of this construction is when $\bO=\DGCat$, $\bo_1=\Vect$, $F_1=\on{Id}$. We will think of $\Phi$ is an 
object $\bc\in \bC$. The datum of $\alpha$ is then a morphism $\bc\to F(\bc)$. 

\smallskip

Assume that $\bc$ is \emph{compact} in $\bC$; this exactly means that the functor $\Vect\to \bC$ defined by $\bc$
admits a right adjoint in $\DGCat$. Thus, to this datum, we can associate a point 
$$\on{cl}(\bc,\alpha)\in \Hom(\Tr(\on{Id},\Vect),\Tr(F,\bC))=\Hom(\sfe,\Tr(F,\bC))=\Tr(F,\bC)).$$

We refer to $\Tr(F,\bC))$ as the \emph{class} of $(\bc,\alpha)$. We use the class construction also more abstractly
for a general $\bO$ with $\bo_1=\one_\bO$ and $F_1=\on{Id}$. 

\ssec{What do we mean by a DG 2-category?} \label{sss:2 cats}

We will not attempt to define what we mean by a general DG 2-category. Rather, we will follow \cite{GRV} and 
consider the \emph{Morita category} of monoidal DG categories. 

I.e., a DG 2-category for us will be always of the form 
$\bA\mod(\DGCat)$, where $\bA$ is a monoidal DG category.
The (2)-category of morphisms $\bA_1\mod(\DGCat)\to \bA_2\mod(\DGCat)$ is by the (2)-category 
$$(\bA_1^{\on{rev-mult}}\otimes \bA_2)\mod(\DGCat).$$

Note that for the above class of DG 2-categories, any object is dualizable. 

\ssec{The 2-categorical class construction} \label{sss:Tr 2 cats} 

Let $\fC:=\bA\mod(\DGCat)$ be a DG 2-category, and let $\fF$ be its endofunctor, given by a bimodule $\bM$. Then 
$$\Tr(\fF,\fC)\simeq \bM\underset{\bA\otimes \bA^{\on{rev-mult}}}\otimes \bA\in \DGCat.$$

Let $\fc\in \fC$ be an object, i.e, $\bC\in \bA\mod(\DGCat)$. A datum of $\alpha$ for such $\fc$ is a map of $\bA$-modules
$$\bC\to \bM\underset{\bA}\otimes \bC.$$

Suppose that $\bC$ is \emph{dualizable} as a $\bA$-module. Then the construction of \secref{sss:class} associates to this
datum an object
$$\on{cl}(\fc,\alpha)\in \Tr(\fF,\fC).$$

\end{document}